\documentclass[aos,MSNbibl,seceqn,nameyear,dvips]{arximspdf}

\doi{10.1214/11-AOS904}
\volume{39}
\issue{5}
\pubyear{2011}
\firstpage{2280}
\lastpage{2301}

\makeatletter

\newcommand{\csi}{\xi}

\renewcommand{\Pi}{\prod}

\newtheorem{prop}{Proposition}
\newtheorem{lem}{Lemma}
\newtheorem{theorem}{Theorem}

\newproclaim{defn}{Definition}

\newcommand{\V}[1]{\mathbf{#1}}
\newcommand{\F}[1]{\mathrm{#1}}

\makeatother

\begin{document}
\begin{frontmatter}

\title{The method of moments and degree distributions for network models}
\runtitle{Method of moments for networks}

\begin{aug}
\author[A]{\fnms{Peter J.} \snm{Bickel}\thanksref{t1}\ead[label=pb]{bickel@stat.berkeley.edu}},
\author[B]{\fnms{Aiyou} \snm{Chen}\thanksref{t2}\ead[label=ac]{aiyouchen@google.com}}
\and
\author[C]{\fnms{Elizaveta} \snm{Levina}\corref{}\thanksref{t3}\ead[label=el]{elevina@umich.edu}}
\runauthor{P. J. Bickel, A. Chen and E. Levina}
\affiliation{University of California, Berkeley,
Google Inc. and University of Michigan}
\dedicated{This research is dedicated to Erich L. Lehmann, the thesis
advisor of one of us and ``grand thesis advisor'' of the others. It is
a work in which we try to develop nonparametric methods for doing
inference in a setting, unlabeled networks, that he never considered.
However, his influence shows in our attempt to formulate and develop a
nonparametric model in this context. We also intend to study to what
extent a potentially ``optimal'' method such as maximum likelihood can
be analyzed and used in this context. In this respect, this is the
first step on a road he always felt was the main one to stick to.}
\address[A]{P. J. Bickel\\
Department of Statistics\\
University of California\\
367 Evans Hall \\
Berkeley, California 94720-3860 \\
USA\\
\printead{pb}}
\address[B]{A. Chen\\
Google Inc.\\
1600 Amphitheatre Pkwy \\
Mountain View, California 94043 \\
USA\\
\printead{ac}}
\address[C]{E. Levina\\
Department of Statistics\\
University of Michigan\\
439 West Hall\\
1085 S. University Ave.\\
Ann Arbor, Michigan 48109-1107 \\
USA\\
\printead{el}}
\end{aug}

\thankstext{t1}{Supported in part by NSF Grant DMS-09-06808.}

\thankstext{t2}{Part of this work was done while A. Chen was at Alcatel-Lucent Bell
Labs.}

\thankstext{t3}{Supported in part by NSF Grants
DMS-08-05798 and DMS-11-06772.}

% HISTORY:
\received{\smonth{2} \syear{2011}}
\revised{\smonth{6} \syear{2011}}

% ABSTRACT
%
\begin{abstract}
Probability models on graphs are becoming increasingly important in
many applications, but statistical tools for fitting such models are
not yet well developed. Here we propose a general method of moments
approach that can be used to fit a large class of probability models
through empirical counts of certain patterns in a graph. We establish
some general asymptotic properties of empirical graph moments and prove
consistency of the estimates as the graph size grows for all ranges of
the average degree including $\Omega(1)$. Additional results are
obtained for the important special case of degree distributions.
\end{abstract}

% KEYWORDS
%

\end{frontmatter}

%s1 ###
\section{Introduction}
\label{secintro}

The analysis of network data has become an important component of doing
research in many fields; examples include social and friendship
networks, food webs, protein interaction and regulatory networks in
genomics, the World Wide web and computer networks. On the algorithmic
side, many algorithms for identifying important network structures such
as communities have been proposed, mainly by computer scientists and
physicists; on the mathematical side, various probability models for
random graphs have been studied. However, there has only been a
limited amount of research on statistical inference for networks, and
on learning the network features by fitting models to data; to a large
extent, this is due to the gap between the relatively simple models
that are analytically tractable and the complex features of real
networks not easily reproduced by these models.

Probability models on infinite graphs have a nice general
representation based on results
[\citet{Aldous81}, \citet{Hoover1979},
\citet{Kallenberg2005}, \citet{Diaconis2008}], analogous to
de Finetti's theorem, for exchangeable matrices. Here, we give a
brief summary closely following the notation of
\citet{BickelChen2009}. Graphs can be represented through their
adjacency matrix $A$, where $A_{ij} = 1$ if there is an edge from node
$i$ to $j$ and 0 otherwise. We assume $A_{ii} = 0$, that is, there are
no self-loops. $A_{ij}$'s can also represent edge weights if the graph
is weighted, and for undirected graphs, which is our focus here,
$A_{ij} = A_{ji}$. For an unlabeled random graph, it is natural to
require its probability distribution $P$ on the set of all matrices
$\{[A_{ij}], i,j \ge1\}$ to satisfy $[A_{\sigma_i \sigma_j}] \sim P$,
where $\sigma$ is an arbitrary permutation of node indices. In that
case, using the characterizations above one can write
%
%e1.1 ###
%
\begin{equation}
A_{ij}
=g(\alpha, \xi_i, \xi_j, \lambda_{ij}) ,
\end{equation}
where $\alpha$,
${\xi_i}$ and $\lambda_{ij}$ are i.i.d. random variables distributed
uniformly on $(0,1)$, $\lambda_{ij}=\lambda_{ji}$ and $g$ is a function
symmetric in its second and third arguments. $\alpha$ as in de
Finetti's theorem corresponds to the mixing distribution and is not
identifiable. The equivalent of the i.i.d. sequences in de Finetti's
theorem here are distributions of the form
$A_{ij} = g(\xi_i, \xi_j, \lambda_{ij})$. This representation is not
unique, and $g$ is not identifiable.
These distributions can be parametrized through the function
%
%e1.2 ###
%
\begin{equation} \label{hclass}
h(u,v) = \mathbb{P}[A_{ij}=1|\xi_i = u, \xi_j = v] .
\end{equation}
The function $h$ is still not unique, but it can be shown that if
two functions $h_1$ and $h_2$ define the same distribution $P$, they
can be related through a measure-preserving transformation, and a
unique canonical $h$ can be defined, with the property that $\int_0^1
h_{\F{can}}(u,v)\,dv$ is monotone nondecreasing in $u$; see
\citet{BickelChen2009} for details. From now on, $h$ will refer to
the canonical $h_{\F{can}}$. We use the following parametrization of
$h$: let
%
%e1.3 ###
%
\begin{equation}
\rho= \int_0^1 \int_0^1 h(u,v) \,du \,dv
\end{equation}
be the
probability of an edge in the network. Then the density of $(\xi_i,
\xi_j)$ conditional on $A_{ij} = 1$ is given by
%
%e1.4 ###
%
\begin{equation}\label{wdef}
w(u,v) = \rho^{-1} h(u,v) .
\end{equation}
With this parametrization, it is natural to
let $\rho= \rho_n$, make $w$ independent of $n$ and control the rate
of the expected degree $\lambda_n = (n-1) \rho_n$ as $n \rightarrow
\infty$.
%(at this point we do not explicitly forbid self-loops, even though
%some commonly studied models do).
The case most studied in probability on random graphs is $\lambda_n =
\Omega(1)$ [where $a_n = \Omega(b_n)$ means $a_n = O(b_n)$ and $b_n =
O(a_n)$]. The case of $\lambda_n = 1$ corresponds to the so-called
phase transition, with the giant connected component emerging for
$\lambda_n > 1$.

Many previously studied probability models for networks fall into this
class. It includes the block model [\citet{Holland83},
\citet{SnijdersNowicki1997}, \citet{Nowicki2001}], the
configuration model
[\citet{ChungLu2002}] and many latent variable models, including the
univariate [\citet{Hoff2002}] and multivariate [\citet
{Handcock2007}]
latent variable models, and latent feature models [\citet
{Hoff2007}]. In
fact, dynamically defined models such as the ``preferential
attachment'' model [which seems to have been first mentioned by Yule in
the 1920s, formally described by \citet{price1965} and given its modern
name by \citet{BarabasiAlbert1999}] can also be thought of in this
way if the dynamical construction process continues forever producing
an infinite graph; see Section 16 of \citet{Bollobas2007}.

\citet{BickelChen2009} pointed out that the block model provides a
natural parametric approximation to the nonparametric model
(\ref{hclass}), and the block model is the main parametric model we
consider in this paper; see more details in Section~\ref{secmodelfit}.
The block model can be defined as follows: each node $i = 1,\ldots, n$
is assigned to one of $K$ blocks independently of the other nodes, with
$\mathbb{P} (c_i = a) = \pi_a$, \mbox{$1\leq a \leq K$}, $\sum_{a=1}^K
\pi_a=1$, where $K$ is known, and ${c} = (c_1,\ldots, c_n)$ is the
$n\times1$ vector of labels representing node assignments to blocks.
Then, conditional on~${c}$, edges are generated independently with
probabilities $\mathbb{P} [A_{ij} = 1|c_i=a, c_j=b] = F_{ab}$. The
vector of probabilities ${\pi} = \{\pi_1,\ldots, \pi_K\}$ and the $K
\times K$ symmetric matrix $F = [F_{ab}]_{1\leq a,b \leq K}$ together
specify a block model. The block model is typically fitted either in
the Bayesian framework through some type of Gibbs sampling
[\citet{SnijdersNowicki1997}] or by maximizing the profile
likelihood using a stochastic search over the node labels
[\citet{BickelChen2009}]. \citet{BickelChen2009} also
established conditions on modularity-type criteria such as the
Newman--Girvan modularity [see \citet{Newman2006} and references
therein] give consistent estimates of the node labels in the block
model, under the condition of the graph degree growing faster than
$\log n$, where $n$ is the number of nodes. They showed that the
profile likelihood criterion satisfies these conditions.%\looseness=-1

The block model is very attractive from the analytical point of view
and useful in a number of applications, but the class (\ref{hclass}) is
much richer than the block model itself. Moreover, the block model
cannot deal with nonuniform edge distributions within blocks, such as
the commonly encountered ``hubs,'' although a modification of the block
model introducing extra node-specific parameters has been recently
proposed by \citet{Karrer10} to address this shortcoming. It may also
be difficult to obtain accurate results from fitting the block model by
maximum likelihood when the graph is sparse.

In this paper, we develop an alternative approach to fitting models of
type (\ref{hclass}), via the classical tool of the method of moments.
By moments, we mean empirical or theoretical frequencies of occurrences
of particular patterns in a graph, such as commonly used triangles and
stars, although the theory is for general patterns. While specific
parametric models like the block model can be fitted by other methods,
the method of moments applies much more generally, and leads to some
general theoretical results on graph moments along the way. We note
that related work on the method of moments was carried out for some
specific parametric models in \citet{Picardetal2008}.

A well-studied class of random graph models where moments play a big
role is the exponential random graph models (ERGMs).
% Note that ERGMs do not belong to our class of exchangeable models;
ERGMs are an exponential family of probability distributions on graphs
of fixed size that use network moments such as number of edges,
$p$-stars and triangles as sufficient statistics. ERGMs were first
proposed by \citet{HollandLeinhardt1981} and \citet
{FrankStrauss1986}
and have then been generalized
in various ways by including nodal covariates or forcing particular
constraints on the parameter space; see \citet{robinsetc2007} and
references therein.
%Unlike the exchangeable models in this paper,
%ERGMs do not specify a probabilistic model for network formation; they
%are merely descriptive, reflecting
%such network properties as transitivity and homophily.
While the ERGMs are relatively tractable, fitting them is difficult
since the partition function can be notoriously hard to estimate.
Moreover, they often fail to provide a good fit to data. Recent
research has shown that a wide range of ERGMs are asymptotically
either too simplistic, that is, they become equivalent to Erd\"
{o}s--Renyi graphs, or nearly degenerate, that is, have no edges or are
complete; see \citet{handcock2003} for empirical studies and
\citet{chatterjeediaconis2011} and \citet
{ShaliziRinaldo2011} for
theoretical analysis.

The rest of the paper is organized as follows. In Section \ref
{secsetup}, we set up the notation and problem formulation and study
the distribution of empirical moments, proving a central limit theorem
for acyclic patterns. We also work out examples for several specific
patterns. In Section \ref{secmodelfit} we show how to use the method of
moments to fit the block model, as well as identify a general
nonparametric model of type (\ref{hclass}).
In Section \ref{secdegree}, we focus on degree distributions, which
characterize (asymptotically) the model (\ref{hclass}). Section \ref
{secbootstrap} discusses the relationship between normalized degrees
and more complicated pattern counts that can be used to simplify
computation of empirical moments. %Section \ref{secnum} gives some
%numerical results, and
Section \ref{secdisc} concludes with a discussion. Proofs and
additional lemmas are given in the \hyperref[app]{Appendix}.
\vfill\eject

%s2 ###
\section{The asymptotic distribution of moments}
\label{secsetup}

%s2.1 ###
\subsection{Notation and theory}

We start by setting up notation. Let $G_{n}$ be a random graph on
vertices $1,\ldots,n$, generated by%\vadjust{\eject}
%
%e2.1 ###
%
\begin{equation}
\label{hcan}
\mathbb{P}(A_{ij}=1|\xi_{i}=u,\xi_{j}=v) = h_{n}(u,v) = \rho_{n}
w(u,v)I(w \leq\rho_n^{-1}),
\end{equation}
where $w(u,v)\geq0$, symmetric, $0\leq u,v\leq1$, $\rho_{n}
\rightarrow0$. We cannot, unfortunately, treat $\rho_n$ and $w$ as two
completely free parameters, as we need to ensure that $h \leq1$. We
can either assume that the sequence $\rho_n$ is such that $\rho_n w
\leq1$ for all $n$, or restrict our attention to classes where
$w_n(u,v) = w(u,v) I(w(u,v) \leq\rho_n^{-1}) \stackrel
{L_2}{\rightarrow
} w(u,v)$. In either case, we can ignore the weak dependence of $w_n$
on $\rho_n$ and effectively replace $w_n$ with $w$.

Let $T\dvtx\mathcal{L}_{2}(0,1)\rightarrow\mathcal{L}_{2}(0,1)$ be the
operator defined by
\[
[Tf](u)\equiv\int_{0}^{1}h(u,v)f(v)\,dv .
\]
We drop the subscript $n$ on $h$, $T$ when convenient. Similarly, let
$T_w\dvtx\mathcal{L}_{2}(0,1)\rightarrow\mathcal{L}_{2}(0,1)$ be
defined by
$w$. Let
\[
D_{i} = \sum_{j}A_{ij},\qquad\bar{D} = \frac{1}{n} \sum
_{i=1}^{n}D_{i} =
\frac{2L}{n} .
\]
Thus $D_i$ is the degree of node $i$, $\bar D$ is the average degree
and $L$ is the total number of edges in $G_n$.

Let $R$ be a subset of $\{(i,j)\dvtx1\leq i<j\leq n\}$. We identify
$R$ with the vertex set $V(R)=\{i\dvtx(i,j) \mbox{ or } (j,i) \in R
\mbox{
for some }j\}$ and the edge set $E(R)=R$. Let $G_n(R)$ be the subgraph
of $G_n$ induced by $V(R)$. Recall that two graphs $R_{1}$ and $R_{2}$
are called isomorphic ($R_{1}\sim R_{2}$) if there exists a one-to-one
map $\sigma$ of $V(R_{1})$ to $V(R_{2})$ such that the map
$(i,j)\rightarrow(\sigma_{i},\sigma_{j})$ is one-to-one from $E(R_{1})$
to $E(R_{2})$.

Throughout the paper, we will be using two key quantities defined next:
\begin{eqnarray*}
Q(R) &=& \mathbb{P}\bigl(A_{ij}=1\mbox{, all }(i,j) \in R\bigr), \\
P(R) &=& \mathbb{P}\bigl(E(G_n(R)) = R\bigr).
\end{eqnarray*}
%
%{\bf[Did not catch this last time, it used to say $P(R) = P(E(G_n) =
%R)$. I think it is correct now.]}
Next, we give a proposition summarizing some simple relationships
between $P$ and $Q$. The proof, which is elementary, is given in the
\hyperref[app]{Appendix}. Similar results are implicit in \citet
{Diaconis2008}.
\begin{prop}
\label{prop1}
If $G_{n}$ is a random graph, and $R$ a subset of $\{(i,j)\dvtx1\leq
i<j\leq n\}$, then
%
%e2.2 ###
%
\begin{eqnarray}\label{eqT2}
P(R) & = & \mathbb{E}\biggl\{ \prod_{(i,j)\in R}h(\xi_{i},\xi
_{j})\prod
_{(i,j)\in\bar{R}}\bigl(1-h(\xi_{i},\xi_{j})\bigr)\biggr\} \nonumber\\ %
& = & Q(R)-\sum\bigl\{Q\bigl(R\cup(i,j)\bigr)\dvtx(i,j)\in\bar{R}\bigr\}\\
& &{} + \sum\bigl\{Q\bigl(R\cup\{(i,j),(k,l)\}\bigr)\dvtx(i,j),
(k,l)\in\bar{R}\bigr\}
- \cdots,\nonumber
% \pm Q(G_{n}(R)).
\end{eqnarray}
where $\bar{R}=\{(i,j)\notin R, i\in V(R), j\in V(R)\}$.
%{\bf[This is still wrong, I think. The second product has to include
%all potential edges between vertices in $V(R)$ that are not in $R$,
%not just those present in $G_n(R)$. It should be ``the complement of
%$R$ in $\{(i,j): i, j \in V(R), \ i < j\}$''.]}
Further,
%
%e2.3 ###
%
\begin{equation}\label{eqT4}
Q(R) = \sum\{P(S)\dvtx S\supset R,V(S)=V(R)\} .
\end{equation}
Here $R\subset S$ refers to $S\subset\{(i,j)\dvtx i,j\in V(R)\}$.
\end{prop}

The quantities $P(R)$ and $Q(R)$ are unknown population quantities
which we can estimate from data, that is, from the graph $G_n$. Define,
for $R\subset\{(i,j)\dvtx1\leq i<j\leq n\}$ with $|V(R)|=p$,
\[
\hat{P}(R) = \frac{1}{{n \choose p}N(R)}\sum\{1(G\sim
R)\dvtx G\subset
G_{n}\},
\]
where $N(R)$ is the number of graphs isomorphic to $R$ on vertices
$1,\ldots, p$. For instance, if $R$ is a 2-star consisting of two edges
$(1,2)$, $(1,3)$, then $N(R) = 3$.
%For instance, if $R$ is a $(k,l)$-wheel (see Definition
Further, let
\[
\hat{Q}(R) = \sum\{\hat{P}(S)\dvtx S\supset R,V(S)=V(R)\}.
\]
Here we use $R$ and $S$ to denote both a subset and a subgraph. Evidently,
\[
\mathbb{E} \hat{P}(R) = P(R) ,\qquad \mathbb{E} \hat{Q}(R) = Q(R).
\]

The scaling here is controlled by the parameter $\rho_n$, the natural
assumption for which is $\rho_{n}\rightarrow0$. In that case,
$P(R)\rightarrow0$
for any fixed $R$ with a fixed number of vertices $p$. Therefore we
consider the following rescaling of $P(R)$ and $Q(R)$: writing $|R|$
for $|E(R)|$, let
\[
\tilde{P}(R) = \rho_{n}^{-|R|}P(R),\qquad \tilde{Q}(R) = \rho
_{n}^{-|R|}Q(R).
\]
Then we have
%
%e2.4 ###
%
\begin{equation}\label{eq4}
\tilde{P}(R) = \mathbb{E} \prod_{(i,j)\in R} w_{n}(\xi_{i},\xi _{j}) +
O\biggl(\frac{\lambda_n}{n}\biggr)
\end{equation}
since
\[
\rho_{n}^{-|R|} \mathbb{ E}\prod_{(i,j)\in R}h_n(\xi_{i},\xi_{j})
\biggl[\prod_{(i,j)\in\bar{R}}\bigl(1-h_{n}(\xi_{i},\xi_{j})\bigr)-1\biggr]
= O(\rho_n) = O\biggl(\frac{\lambda_n}{n}\biggr),
\]
if $\int w^{2(|R|+1)}(u,v)\,du \,dv<\infty$.

Next, we define the natural sample estimates of the population
quantities $\tilde P$ and $\tilde Q$ by
\[
\check{P}(R) = \hat\rho_n^{-|R|}\hat{P}(R),\qquad
\check{Q}(R) = \hat\rho_n^{-|R|}\hat{Q}(R) ,
\]
where\vspace*{1pt} $\hat\rho_n = \frac{\bar{D}}{n-1} = \frac{2L}{n(n-1)}$ is the
estimated probability of an edge. For these rescaled versions of $P$
and $Q$, we have the following theorem.
\begin{theorem}
\label{thmLandP} Suppose $\int_{0}^{1}\int_{0}^{1}w^{2}(u,v)\,dv
\,du<\infty$.

\begin{longlist}[(a)]
\item[(a)]
If $\lambda_n\rightarrow\infty$, then
%
%e2.6 ###
%e2.5 ###
%
\begin{eqnarray}
\label{eqneweq5}
\frac{\hat\rho_n}{\rho_n} &\rightarrow_{P}& 1, \\
\label{eqneweq6}
\sqrt{n}\biggl(\frac{\hat\rho_n}{\rho_n} - 1 \biggr)&\Rightarrow&
\mathcal{N} (0,\sigma^{2})
\end{eqnarray}
for some $\sigma^2 > 0$.
Suppose further $R$ is fixed, acyclic with $|V(R)|=p$ and $\int
w^{2|R|}(u,v)\,du \,dv<\infty$. Then,
%
%e2.7 ###
%
\begin{eqnarray}\label{eqneweq7}
\check{P}(R)&\to_{P}&\tilde{P}(R) , \nonumber\\[-8pt]\\[-8pt]
\sqrt{n}\bigl(\check{P}(R)-\tilde{P}(R)\bigr) &\Rightarrow&
\mathcal{N}(0,\sigma
^{2}(R)) .\nonumber
\end{eqnarray}
More generally, for any fixed $\{R_{1},\ldots,R_{k}\}$ as above with
$|V(R_{j})|\leq p$,
%
%e2.8 ###
%
\begin{equation}\label{eqneweq8}
\sqrt{n}\bigl((\check{P}(R_{1}),\ldots,\check{P}(R_{k}))-(\tilde
{P}(R_{1}),\ldots,\tilde{P}(R_{k}))\bigr)
\Rightarrow\mathcal{N}(\V{0},\Sigma(\V{R})) .
\end{equation}
\item[(b)] Suppose $\lambda_n \rightarrow\lambda < \infty$. Conclusions
(\ref{eqneweq5})--(\ref{eqneweq8}) continue to hold save that
$\sigma^{2}(R)$, $\Sigma(R)$ depend on $\lambda$ as well as $R$.

\item[(c)] Even if $R$ is not necessarily acyclic, the same conclusions
apply to $\check{Q}$ and $\tilde{Q}$ if $\lambda_n$ is of order
$n^{1-2/p}$ or higher, and to $\check{P}$ and $\tilde{P}$ under the
same condition on~$\lambda_n$.
\end{longlist}
\end{theorem}

The proof is given in the \hyperref[app]{Appendix}.

\subsubsection*{Remarks}

(1) Note that part (b) yields consistency and asymptotic normality of
acyclic graph moment
estimates across the phase transition to a giant component, that is,
for $\lambda<1$ as well
as $\lambda\geq1$.

(2) Note that we are, throughout, estimating features of the canonical $w$.
Unnormalized $P$ and $Q$ are trivially 0 if $\lambda_n$ is not of
order $n$.

(3) In view of (\ref{eq4}), we can use $\check{P}(R)$ as an
estimate of $\tilde{Q}(R)$ if $R$ is acyclic and $\lambda
_n=o(n^{1/2})$, since in this case the bias of $\check{P}$ is of order
$o(n^{-1/2})$. The reason for not using $\check{Q}(R)$ directly even if
$R$ is acyclic is that by (\ref{eqT4}),
there may exist $S\supset R$ which are not acyclic, and we can
therefore not
conclude that the theorem also applies to $\check{Q}$ unless we are in
case (c).

%Another difficulty of using $\check{Q}$ is that, except for some
%special patterns discussed below, there is no nice general formula for
%computing $\check{P}(S)$, and we need to evaluate $\check{Q}(R)$ using

(4) Part (c) of the theorem shows that for graphs
with $\lambda_n=\Omega(n)$, $\check{Q}$ always gives $\sqrt{n}$-consistent
estimates of any pattern while $\check{P}$ is not consistent unless we
assume acyclic graphs,\vspace*{1pt} since the bias is of order $O(\lambda_n/n) =
O(1)$. In the
range $\lambda_n=o(n^{1/2})$ to $\Omega(n)$, what is possible
depends on the pattern. For instance, if $\Delta=\{(1,2),(2,3),(3,1)\}$,
a triangle, $\check{P}(\Delta)=\check{Q}(\Delta)$ (because there is
no other
graph on three nodes containing $\Delta$), and $\check{P}$ is
$\sqrt{n}$-consistent if $\lambda_n\geq\varepsilon n^{1/3}$ by part (c) but
otherwise only consistent if $\lambda_n\to\infty$.

%s2.2 ###
\subsection{Examples of specific patterns}
Next we give explicit formulas for several specific $R$. Our main focus
is on wheels (defined next), which, as we shall see, in principle can
determine the canonical $w$.
%give particularly concise moment formulas.
%
\begin{defn}[(Wheels)]
\label{defkl-wheel}A $(k,l)$-wheel is a graph with $kl+1$ vertices
and\break $kl$ edges isomorphic to the graph with edges $\{(1,2),\ldots
,(k,k+1); (1,k+2),\break\ldots,(2k,2k+1);\ldots,(1,(l-1)k+2),\ldots
,(lk,lk+1)\}$.
\end{defn}

In other words, a wheel consists of node $1$ at the center and
$l$ ``spokes'' connected to the center, and each spoke is a chain of
$k$ edges. We consider only $k\geq2$. The number of isomorphic
$(k,l)$-wheels on vertices $1,\ldots, p$ is $N(R)= (kl+1)!/l!$.

If the graph $R$ is a $(k,l)$-wheel, the theoretical moments have a
simple form and can be expressed in terms of the operator $T$ as follows:
%
%e2.9 ###
%
\begin{equation}\label{prokl-wheel}
Q(R) = \mathbb{E} (T^{k}(1)(\xi_{1}))^{l} .
\end{equation}
This follows from
\begin{eqnarray*}
Q(R) & = & \mathbb{E}\Bigl(\mathbb{E}\Bigl(\prod\{h(\xi_{i},\xi
_{j})\dvtx(i,j)\in
E(R)\}|\xi_{1}\Bigr)\Bigr) \\
& = & \biggl(\int_{0}^{1}\cdots\int_{0}^{1}h(\xi_{1},\xi
_{2})\cdots
h(\xi_{k},\xi_{k+1})\,d\xi_{2}\cdots d\xi_{k+1}\biggr)^{l} \\
& = & \mathbb{E} (T^{k}(1)(\xi_{1}))^{l},
\end{eqnarray*}
where the first equality holds by the definition of $Q$ and the second
by the structure of a $(k,l)$-wheel.

%From our general considerations, consistency in the sense of
%convergence to $\tilde Q(R)$ always holds for $\check P$ {\bf[should
%this be $\check Q$? ]} applied to $(k,l)$-wheels and $\lambda_n=o(n)$,
%but $\sqrt{n}$-consistency
%appears to hold in general only if $\lambda_n=o(n^{1/2})$. In general,
%$\check Q$ is $\sqrt{n}$-consistent for $\lambda_n=\Omega(n)$.
%Otherwise it would appear that $\check{P}$
%has bias of order $\geq cn^{-1/2}$ while $\check{Q}$
%has variance of order larger than $n^{-1}$ unless $\lambda_n\geq
%2/(kl+1)}$. For example, let $k=l=2$. Then for $\check P(R)$ to be a
%consistent estimate of $\tilde Q(R)$ we need to have $\rho_n E(\check
%P(\Delta^2) = o(n^{-1/2})$ even if $\lambda_n = \Omega(n^{1/2})$.
%Suppose
%$\Delta^2$ is a figure, two triangles with a common edge. This is the
%order of $E\check P(R) - E\check Q(R)$. But $E \check P(\Delta^2) =
%O(1)$ and $\lambda_n = \Omega(n^{1/2})$ means that $\rho_n =
%Is this a typo and should it be $\frac{n^{-4}}{\rho_n^{10}}$? Still
%can't follow it. This paragraph needs to be rewritten more clearly and
%perhaps more formally, Peter, could you please give it a try? ]}

For a $(k,l)$-wheel $R$, from our general considerations, $\mathbb{E}
\check{P}(R) = \tilde{P} (R) = \tilde{Q}(R) + o(1)$ if $\lambda_n=o(n)$
and in view of (\ref{eqneweq8}), $\check{P}(R)$ always consistently
estimates $\tilde{Q}(R)$. However, $\sqrt{n}$-consistency of $\check
{P}$ (converging to $\tilde{Q}$) holds in general only if $\lambda_n =
o(n^{1/2})$. By part (c) $\check{Q}$ is $\sqrt{n}$ consistent for
$\tilde{Q}$ only if $\lambda_n$ is of\vspace*{1pt} order larger than
$n^{1-2/(kl+1)}$. In the $\lambda_n$ range between $O(n^{1/2})$ and
$O(n^{1-2/(kl+1)})$, we do not exhibit a $\sqrt{n}$-consistent estimate
though we conjecture that by appropriate de-biasing of $\check{P}$ such
an estimate may be constructed. However, $\lambda_n=o(n^{1/2})$ seems a
reasonable assumption\vspace*{1pt} for most graphs in practice, and
then we can use the more easily computed $\check{P}$.
\begin{defn}[(Generalized wheels)]
\label{crosswheel} A $(\V{k},\V{l})$-wheel, where $\V
{k}=(k_{1},\ldots,k_{t})$,
$\V{l}=(l_{1},\ldots,l_{t})$ are vectors
and the $k_{j}$'s are distinct integers, is the union $R_1 \cup\cdots
\cup R_t$, where $R_j$ is a $(k_j, l_j)$-wheel,
$j = 1,\ldots, t$, and the wheels $R_1,\ldots, R_t$ share a common hub
but all their spokes are disjoint.
\end{defn}

A\vspace*{1pt} $(\V{k},\V{l})$-wheel has a total of $p=\sum_{j}l_{j}k_{j}+1$
vertices and $\sum_{j}l_{j}k_{j}$ edges. For example, a graph defined by
$E=\{(1,2);(1,3),(3,4);(1,5),(5,6);(1,7),(7,8)$, $(8,9)\}$
is a $(\V{k},\V{l})$-wheel with $\V{k}=(1,2,3)$ and
$\V{l}=(1,2,1)$. The number of distinct isomorphic $(\V{k},\V{l})$-wheels
on $p$ vertices is $N(R) = p! (\prod_{j}l_{j}!)^{-1}$.\vspace*{1pt}

We can compute, defining $A(R)=\Pi\{A_{ij}\dvtx(i,j)\in R\}$,
%
%e2.10 ###
%
\begin{eqnarray}\label{eq15}
Q(R) & = & \mathbb{P}\Biggl(\bigcap_{j=1}^{t}[A(R_{j})=1]
\Biggr)\nonumber\\
& = & \mathbb{E}\Biggl\{\prod_{j=1}^{t}\mathbb{P}\bigl(A(R_{j})=1|
\mbox{Hub}\bigr)\Biggr\} \\
& = & \mathbb{E}\prod_{j=1}^{t}[T^{k_{j}}(\xi
)]^{l_{j}}.\nonumber
\end{eqnarray}
Thus $(\V{k},\V{l})$-wheels give us all cross moments of
$T^{m}(\xi)$, $m \geq1$. Note that all $(\V{k},\V{l})$-wheels are acyclic.

We are not aware of other patterns for which the moment formulas are as
simple as those for wheels. For example, if $R$ is a triangle, then
\begin{eqnarray*}
Q(R) & = & \int_{0}^{1}\int_{0}^{1}\int_{0}^{1}h(u,v)h(v,w)h(w,u)\,du
\,dv
\,dw\\
& = & \int_{0}^{1}\int_{0}^{1}h^{(2)}(u,w) h(w,u) \,du \,dw,
\end{eqnarray*}
where $h^{(2)}(u,w) = \int_0^1 h(u,v) h(v,w) \,dv$ corresponds to
$T^{2}f\equiv\int_{0}^{1}h^{(2)}(u,v)\times\break f(v) \,dv$.

In general, unions of $(\V{k},\V{l})$-wheels are also more complicated.
If $R_{1},R_{2}$ are $(\V{k}_{1},\V{l}_{1})$,
$(\V{k}_{2},\V{l}_{2})$-wheels which share a single node
[$V(R_{1})\cap
V(R_{2})=\{a\}$], we can compute $P(R_{1}\cup R_{2})=\mathbb{E}
P(R_{1}|\xi_{a})P(R_{2}|\xi_{a})$. If $a$ is the hub of both wheels,
then evidently $R_{1}\cup R_{2}$ is itself a generalized wheel, and
(\ref{eq15}) applies. Otherwise, the formula, as for triangles, is more
complex. However, such unions of $(\V{k},\V{l})$-wheels are acyclic.

%s3 ###
\section{Moments and model identifiability}
\label{secmodelfit}

We establish two results in this section: identifiability of block
models with known $K$ using $\{\check{P}(R)\dvtx R$ a $(k,l)$-wheel,
$1\leq l\leq2K-1,2\leq k\leq K\}$, and\vspace*{1pt} the general
identifiability of the function $w$ from $\{\check{P}(R)\}$ using all
$(\V{k},\V{l})$-wheels $R$.

%s3.1 ###
\subsection{The block model}

Let $w$ correspond to a $K$-block model defined by parameters $\theta
\equiv(\pi,\rho_{n},S)$, where $\pi_{a}$ is the probability of a node
being assigned to block $a$ as before, and
\[
F_{ab} \equiv\mathbb{P} (A_{ij}=1|i\in a,j\in b)= \rho_{n}S_{ab},
\qquad 1\leq a,b\leq K .
\]

Recall that the function $h$ in (\ref{hclass}) is not unique, but a
canonical $h$ can be defined. For the block model, we use the
canonical $h$ given by \citet{BickelChen2009}. Let $H_{ab} =
S_{ab}\pi_a\pi_b$. Let the labeling of the communities $1,\ldots, K$
satisfy $H_1 \leq\cdots\leq H_K$, where $H_a = \sum_b H_{ab} $ is
proportional to the expected degree for a member of block $a$. The
canonical function $h$ then takes the value $F_{ab}$ on the $(a,b)$
block of the product partition where each axis is divided into
intervals of lengths $\pi_1,\ldots, \pi_K$. Let $F \equiv
\|F_{ab}\|$.

In view of (\ref{eqneweq6}), we will treat $\rho_{n}$ as known. Let
$\{
W_{kl}\dvtx1\le l\le2K-1,2\leq k\leq K\}$ be the specified set
of $(k,l)$-wheels, and let
\[
\tau_{kl}=\rho^{-kl}P(W_{kl})=\tilde{P}(W_{kl}) ,\qquad \check{\tau
}_{kl}=\check{P}(W_{kl}).
\]
%
%{\bf[should the scaling be $\rho^{-kl-1}$, since a $(k,l)$ wheel has
%$kl+1$ nodes?]}
Let $f\dvtx\Theta\to\mathbb{R}^{(2K-1)(K-1)}$ be the map carrying
the parameters of the block model $\theta\equiv(\pi,S)$ to $\tau
\equiv
\|\tau_{kl}\|$.
$\Theta$ here is the appropriate open subset of $\mathbb{R}^{K(K+3)/2
-2}$. Note that the number of free parameters in the block model is
$K-1$ for $\pi$ and
$K(K+1)/2$ for $F$, but $S$ only has $K(K+1)/2 -1$ free parameters, to
account for $\rho$.
\begin{theorem}
\label{thm2}
Suppose $\theta=(\pi,S)$ defines a block model with known $K$, and
the vectors
$\pi,F\pi,\ldots,F^{K-1}\pi$ are linearly independent. Suppose
$\varepsilon
\leq\lambda_n=o(n^{1/2})$. Then:

\begin{longlist}[(a)]
\item[(a)] $\{\tau_{kl}\dvtx l=1,\ldots,2K-1,k=2,\ldots,K\}$ identify
the $K(K+3)/2-2$ parameters of the block model other than $\rho$ (i.e.,
the map $f$ is one to one).\vspace*{1pt}

\item[(b)] If $f$ has a gradient which is of rank $\frac{K(K+3)}{2}-2$
at the true $(\pi_{0},S_{0})$, then ${f}^{-1}(P(\check{\tau}))$ is a
$\sqrt {n}$-consistent estimate of $(\pi_{0},S_{0})$, where
$\check{\tau}=\|\check{\tau }_{kl}\| $ and $P(\check{\tau})$ is the
closest point in the range of $f$ to $\check{\tau}$.
\end{longlist}
\end{theorem}

Note that the linear independence condition rules out all matrices $F$
that have~$1$ as an eigenvector.  In particular, it rules out the case
of $F_{aa}$ equal for all $a$, $F_{ab}$ equal for all $a \neq b$, which
was studied in detail by \citet{Decelleetal2011}. Using physics
arguments, they showed that in that particular case, when $\lambda =
O(1)$, there are regions of the parameter space where neither the
parameters nor the block assignments can be estimated by any method.

Part (b) shows $\sqrt{n}$-consistency of nonlinear least squares estimation
of $(\pi,S)$ using $\check{\tau}$ to estimate $\tilde{\tau}(\theta,S)$.
The variance of $\check{\tau}_{kl}$ is proportional asymptotically
to that of $\mathbb{E}\{\prod_{(i,j)\in S}w(\xi_{i},\xi_{j})|\xi
_{1}\}$,
where $\xi_{1}$ corresponds to the hub, which we expect increases
exponentially in $p=kl+1$.
If we knew these variances, we could use weighted nonlinear
least squares. In Section \ref{secbootstrap}, we suggest a bootstrap method
by which such variances can be estimated, but we do not pursue this
further in this paper.

%s3.2 ###
\subsection{The nonparametric model}

In the general case, we express everything in terms of the operator $T_{w}
\equiv T/\rho_n$ induced by the canonical $w$. We require that:

\begin{longlist}[(A)]
\item[(A)]
the joint distribution of $\{T_w^{l}(1)(\xi)\dvtx l\geq1\}$
is determined by the cross moments of $(T_w^{l_{1}}(\xi),\ldots
,T_w^{l_{k}}(\xi))$, for $l_{1},\ldots,l_{k}$ arbitrary.
\end{longlist}

A simple sufficient condition for (A) is $|w|\leq M<\infty$. A more
elaborate one is the following:
\begin{longlist}[(A$'$)]
\item[(A$'$)]
\begin{eqnarray*}
\\[-26pt]
\mathbb{E} e^{sw^{k}(\xi_{1},\xi_{2})}&<&\infty,\qquad
0\leq
|s|\leq\varepsilon\mbox{ all }k\mbox{ some }\varepsilon>0 .
\end{eqnarray*}
\end{longlist}
\begin{prop}
\label{propA}
Condition \textup{(A$'$)} implies \textup{(A)}.
\end{prop}

The proof is given in the \hyperref[app]{Appendix}.

Let $w$ characterize $T_w$, where $\int_{0}^{1}w^{2}(u,v)\,du
\,dv<\infty$.
By Mercer's theorem,
%
%e3.1 ###
%
\begin{equation}
\label{eigenw}
w(u,v)=\sum_{j}\lambda_{j}\phi_{j}(u)\phi_{j}(v),
\end{equation}
where the $\phi_{j}$ are orthonormal eigenfunctions and the $\lambda_{j}$
eigenvalues, \mbox{$\sum\lambda_{j}^{2}<\infty$.}\looseness=-1
\begin{theorem}
\label{thmnonparametric} Suppose $\int_{0}^{1}\int
_{0}^{1}w^{2}(u,v)\,du \,dv<\infty$.
Assume the eigenvalues $\lambda_{1}>\lambda_{2}>\cdots$
of $T_w$
are each of multiplicity $1$ with corresponding eigenfunction
$\phi_{j}$, and $\int_{0}^{1}\phi_{j}(u)\,du\neq0$
for all $j$. The joint distribution of $(T_w(1)(\xi),\ldots,\break
T_w^{m}(1)(\xi),\ldots)$
then determines, and is determined by, $w(\cdot,\cdot)$.
\end{theorem}

Note again that interesting cases are ruled out by the condition that
all eigenfunctions of $T$ are not orthogonal to $1$. The general
analogue to the  block model case is that $P(A_{ij}=1|\csi_i)$ cannot
be constant for all $i$ and $j$. Constancy can be interpreted as saying
that $A_{ij}$ and the latent variable $\csi_i$ associated with vertex
$i$ are independent. The proof of Theorem \ref{thmnonparametric} is
given in the \hyperref[app]{Appendix}. The almost immediate application
to wheels is stated next.

%The proof is given in the \hyperref[app]{Appendix}. The almost
%immediate application to
%wheels is stated next.
%
\begin{theorem}
\label{thmnonparaconsistency} Suppose assumption \textup{(A)} and the conditions
of Theorem \ref{thmnonparametric} hold. Let
$\tau_{\V{kl}}=\tilde{P}(S_{\V{kl}})$ where $S_{\V{kl}}$ is a $(\V{k},
\V{l})$-wheel. Then $\mathcal{S}\equiv\{\tau_{\V{kl}}\mbox{: all
}\V{k},
\V{l}\}$ determines $T$.
If $\check{\tau}_{\V{kl}}\equiv\check{P}(S_{\V{kl}})$,
$\check{\tau}_{\V{kl}}$ are $\sqrt{n}$-consistent estimates
of $\tau_{\V{kl}}$, provided that $\lambda_n=o(n^{1/2})$.
\end{theorem}
\begin{pf}%[Proof of Theorem \ref{thmnonparaconsistency}]
Since $\V{T}_{l}\equiv(T(1)(\xi),\ldots,T^{l}(\xi))$
has a moment generating function converging on $0<|s|\leq\varepsilon_{l}$,
the moments (including cross moments) determine the distribution of
the vector.
% by e.g. Theorems in Feller II \citet{} {\bf needs to add a
%citation?}.
%Theorem \ref{thmnonparametric}.
By (\ref{eq15}), the $\tau_{\V{kl}}$ give all moments
of the vector $\V{T}_{l}$ for all $l$. By Theorem \ref{thmLandP}, the
$\check{\tau}_{\V{kl}}$ are $\sqrt{n}$-consistent.
\end{pf}

%s4 ###
\section{Degree distributions}
\label{secdegree}

The average degree $\bar{D}$ is, as we have seen in Theorem
\ref{thmLandP}, a natural data dependent normalizer for moment\vadjust{\eject}
statistics which eliminates the need to ``know'' $\rho_n$. In fact, as we
show in this section, the joint empirical distribution of degrees and
what we shall call $m$ degrees below can be used in estimating
asymptotic approximations to $w(\cdot,\cdot)$ in a somewhat more direct
way than moment statistics. They can also be used to approximate moment
estimates based on $(\V{k}, \V{l})$-wheels in a way that potentially
simplifies computation.

We define\vspace*{-1pt} the $m$-degree of $i$, $D_{i}^{(m)}$, as
the total number of loopless paths of length $m$ between $i$ and
other vertices. Note that the $D_{i}^{(m)}$ can be interpreted as
the ``volume'' of the radius $m$ geodesic sphere around $i$. As for
regular degrees, we normalize and consider
$D_{i}^{(m)}/\bar{D}^{m}$, $i=1,\ldots,n$, and
the empirical joint distribution of vectors $\V{D}_{i}^{(m)}\equiv
(\frac{D_{i}}{\bar{D}},\frac{D_{i}^{(2)}}{\bar{D}^{2}},\ldots
,\frac
{D_{i}^{(m)}}{\bar{D}^{m}})$,
$i=1,\ldots,n$.
The generalized degrees can be computed as follows: for all entries of $A^m$,
eliminate all terms in the sum defining each entry in which an
index appears more than once to obtain a modified matrix $\tilde
{A}^{(m)}=[\tilde{A}_{ij}^{(m)}]$; then
the $D_i^{(m)}$ are given by row sums of $\tilde{A}^{(m)}$. In other
words, letting $A_{E(R)} = \prod_{(i,j)\in E(R)}A_{ij}$ we can write
\begin{eqnarray*}
\tilde{A}_{ij}^{(m)} &=& \sum\bigl\{A_{E(R)}\dvtx R=\{
(i,i_{1}),(i_{1},i_{2}),\ldots,(i_{m-1},j)\}, \\
&&\hspace*{100.2pt} i, i_{1},\ldots,i_{m-1}, j \mbox{ distinct} \bigr\} .
\end{eqnarray*}

The complexity of this computation
% the same as that of computing efficiently $\{D_i^{(m)}: 1\leq i\leq n
is $O((n+m)\lambda_n^m)$ (first term is for computing the row sums of
$A^m$ and the second for eliminating the loops).\vadjust{\goodbreak}

Define the empirical distribution of the vector of normalized degrees
\[
\hat{F}_{m}(\V{x}) = \frac{1}{n}\sum_{i=1}^{n}1\bigl(\V
{D}_{i}^{(m)}\leq\V
{x}\bigr).
\]
Further,\vspace*{1pt} recall the Mallows 2-distance between two distributions $P$
and $Q$, defined by $M_2(P,Q) = \min_F \{ (\mathbb{E}
\|X-Y\|^2)^{1/2}\dvtx
(X, Y) \sim F, X \sim P, Y \sim Q \}$. A~sequence of distribution
functions $F_n$ converges to $F$ in $M_2$ ($F_n \stackrel
{M_2}{\rightarrow} F$) if and only if $F_n \Rightarrow F$ in
distribution, and $F_n$, $F$ have second moments such that
$\int|\V{x}|^2 \,dF_n(\V{x}) \rightarrow\int|\V{x}|^2 \,dF(\V{x})$.
\begin{theorem}
\label{thmldegree} Suppose $\lambda_n\to\infty$ and
$|w_{2m}|<\infty$.
Then ${\hat{F}}_{m} \stackrel{M_2}{\rightarrow} F_{m}$ as $n
\rightarrow\infty$, where $F_{m}$ is the distribution of
$\bolds{\theta}_{m}(\xi)=(\tau_w(\xi),\ldots,T_w^{m-1}(\tau
_w)(\xi))$,
and $\tau_w(\xi)=\int_{0}^{1}w(\xi,v)\,dv$ is monotone increasing.
Moreover, if $\hat G_m(\V{x}, \V{y})$ is the empirical distribution of
$(\V{D}_i^{(m)}, \bolds{\theta}_m(\xi_i))$, then
%
%e4.1 ###
%
\begin{equation}\label{eqmallowsjoint}
\int| \V{x} - \V{y} |^2 \,d \hat G_m(\V{x}, \V{y}) \stackrel
{P}{\rightarrow} 0 .
\end{equation}
\end{theorem}

The proof is given in the \hyperref[app]{Appendix}.
\eject

There is an attractive interpretation of the last statement of Theorem
\ref{thmldegree}.
If $\lambda_n\to\infty$, $\lambda_n=o(n^{1/(m-1)})$, $m\geq2$,
then $D_{i}/\lambda_n$ can be identified with $\tau(\xi_{i})$ in the
following sense:
While $\xi_i$ is unobserved but $D_i/\bar{D}$ is, on average,
$\tau
(\xi_i)$ and $D_i/\bar{D}$ are close. Since $\tau$ is monotone
increasing in $\xi$, that is, is a measure of $\xi$ on another scale,
% {\bf[I still don't understand what the phrase ``can be identified''
%means in this context.]}
we can treat $D_{i}/\lambda_n$ as the latent affinity
of $i$ to form relationships.

\citet{Bollobas2007} show that if
$m=1$, $\lambda_n=O(1)$, then the limit of the empirical distribution
of the degrees can be described as follows: given $\xi\sim\mathcal
{U}(0,1)$, the limit distribution is Poisson with mean $\tau_w(\xi)$.
The limit of the joint degree distribution in this case can
be determined but does not seem to give much insight.

\subsection*{Remark}

Theorem \ref{thmldegree} shows that the normalized degree
distributions can be used for estimation of parameters only if $\lambda
_n\to\infty$. If that
is the case we can proceed as follows:

\begin{longlist}[(1)]
\item[(1)] Let $\hat{\tau}_{1},\ldots,\hat{\tau}_{n}$ be the empirical
quantiles of the normalized 1-degree distribution, and let
$\hat{T}^{m}(\hat {\tau }_{k})$ be the $m$-degree of the vertex with
normalized degree~$\hat {\tau}_{k}$.
\item[(2)] Fit smooth curves to $(\hat{\tau}_{k},\hat
{T}^{m}(\hat
{\tau}_{k}))$ viewed as observations of functions at $\hat{\tau
}_{k}$, $k=1,\ldots,n$, for each $m$, and call these $\hat
{T}^{m}(\cdot
)$ (on $R$).
By Theorem \ref{thmldegree}, $\hat{T}^{m}(t)\rightarrow T^{m-1}(\tau
)(\tau^{-1}(t))$
for all $t$. If $T^{m-1}(\tau^{-1}(\cdot))$ are smooth, the
convergence can
be made uniform on compacts.
\item[(3)] From the fitted functions $\hat{T}^{m}(\cdot)$, we can
estimate the
parameters of block models of any order consistently by replacing
$\V{v}_{m}$ in the proof of identifiability of block models by fitting the
$\hat{T}^{m}(t)$ by $T^{m}(t)$ of the
type specified by block models and then using the corresponding $\hat
{\V
{v}}_{m}$. We only need the conditions of Theorem \ref{thmldegree}.
\end{longlist}

\section{Computation of moment estimates and estimation of their variances}
\label{secbootstrap}

General acyclic graph moment
estimates including those corresponding to patterns arising from
$(\V{k},\V{l})$-wheels are computationally difficult. For
$(k,l)$-wheels with small $k$ and $l$, we can use brute force counting,
but unfortunately,
the complexity of moment computation even for $(k,l)$-wheels
appears to be $O(n\lambda_n^{k})$. Note that we need to count the
sets of loopless paths of length $k$, $S_{i\V{a}}$, for each $i$, where
$S_{i\V{a}}$ is the set of all paths of length $k$ originating at node
$i$ which intersect another such path at $a_1 < \cdots< a_m$, $1 \leq m
\leq k$, and $S_{i0}$ is the set of all paths of length $k$ from $i$
which do not intersect. The number
of $(k,l)$-wheels with hub $i$ is then the number of $l$-tuples
of such paths selected so that elements from $S_{i\V{a}}$ appear at
most once, with the remaining paths coming from $S_{i0}$. This is
computationally nontrivial.

For very sparse graphs, however, intersecting paths can be ignored up
to a certain order, and the wheel counts can be related to normalized
$m$-degrees via a following approximation. If the conditions of Theorem
\ref{thmldegree} hold and $\lambda_n=o(n^{\alpha})$ for all $\alpha
>0$, then
%
%e5.1 ###
%
\begin{equation}\label{eqinthm5}
\hat{\tau}_{kl} = \frac{1}{n}\sum_{i=1}^{n}\frac
{(D_{i}^{(k)})_{l}}{\bar{D}^{kl}}+o_{P}(n^{-1/2}).
\end{equation}
%
%where $(a)_{b}$ is the number of permutations of $b$ things out of
%$a$, $(a)_b = a(a-1)\cdots(a-b+1)$.
A similar formula holds for $\hat{\tau}_{\V{k}\V{l}}$.

The heuristic argument for (\ref{eqinthm5}) is that the expected number
of paths of lengths $k$ from $i$ is $O(\lambda_n^k)$. The expected
number of pairs of such paths which intersect at least once is
\begin{eqnarray*}
&& O(\lambda_n^{2k}) \mathbb{P} [\mbox{two specified paths intersect at
least once}] \\
&&\qquad= O\bigl(\lambda_n^{2k}\bigl( 1 - (1 - \lambda_n/n)^k\bigr)\bigr) =
O\biggl( \frac{k\lambda_n^{2k+1}}{n} \biggr) = o(1),
\end{eqnarray*}
if $\lambda_n = o(n^{\alpha})$ for all $\alpha> 0$. Note that for
$K$-block models this condition is not necessary for all $\alpha$,
since we only need to count a finite number of $(k,l)$-wheels.

%If $\lambda_n=o(n^{1/4})$ the equivalence in \eqref{eqinthm5} holds
%to order $o_{P}(n^{-1/2})$ and hence we can use the formula
%for inference at the $n^{-1/2}$ scale.

Estimation of variances of moment estimates even for $(\V{k},\V{l})$-wheels involve the counting of more complicated patterns. However,
we propose the following bootstrap method:

\begin{longlist}
\item
Associate with each vertex $i$ the counts of $(\V{k}, \V{l})$-wheels
for which it is a hub, $S_{i}=\{n_{i\V{k}\V{l}}\mbox{: all } \V{k},
\V
{l}\}$,
$i=1,\ldots,n$.

\item Sample without replacement $m$ vertices $\{i_{1},\ldots,i_{m}\}$,
and let
\[
\bar{D}^{*} = \frac{1}{m}\sum_{j=1}^{m}D_{i_{j}}.
\]
%
%Let $G_{m}^{*}$ be the subgraph of $G_n$ induced by the selected
%vertices.
%and the associated sets $\{S_{i_{j}}\}$, $j=1,\ldots,m$.
For $R$ a $(\V{k},\V{l})$-wheel, define
\begin{eqnarray*}
\hat{P}^{*}(R) & = & \frac{({n/m})\sum_{j=1}^{m}n_{i_{j}\V
{k}\V
{l}}}{{n \choose p}N(R)}, \\
\check{P}^{*}(R) & = & \hat{P}^{*}(R)\biggl(\frac{\bar
{D}^{*}}{m}\biggr)^{-|R|}.
\end{eqnarray*}

\item Repeat this $B$ times to obtain $\check{P}_{1}^{*},\ldots
,\check
{P}_{B}^{*}$, and let
\[
\hat{\sigma}^{2} = \frac{m}{n}\frac{1}{B}\sum_{b=1}^{B}(\check
{P}_{b}^{*}-\check{P}_{\cdot}^{*})^{2}.
\]
\end{longlist}
Then $\hat{\sigma}^{2}$ is an estimate of the variance of $ \check{P}(R)$
if $\frac{m}{n}\rightarrow0,m\rightarrow\infty$.

This scheme\vspace*{1pt} works if $\lambda_n\rightarrow\infty$ since, given that
the first term of $\check{P}(R)-\tilde{P}(R)$ is of lower order given
$\xi_{1},\ldots,\xi_{n}$, each $\tilde{P}^{*}(R)$ corresponds to
a sample without replacement from the set of possible $\{\xi_{i}\}$.
We conjecture that this bootstrap still works if $\lambda_n=O(1)$.
A similar device can be applied to approximation (\ref{eqinthm5}).

%s6 ###
\section{Discussion}
\label{secdisc}

%s6.1 ###
\subsection{Estimation of canonical $w$ generally}
Our Theorem \ref{thmnonparaconsistency} suggests that we
%could, by letting the number of blocks in block models grow, be able to
might be able to construct consistent nonparametric estimates of
$w_{\mathit{CAN}}$. That is, $\bolds{\tau}_M = \{\tau_{\V{k}\V{l}} \dvtx|\V{k}|
\leq M, |\V{l}| \leq M \}$ can be estimated at rate $n^{-1/2}$ for all
$M < \infty$. But $\{ \bolds{\tau}_M, M \geq1\}$ determines $T_w$, and
thus in principle we can estimate $T_w$ arbitrarily closely using $\{
\hat\tau_{\V{k}\V{l}} \}$. This appears difficult both theoretically
and practically. Theoretically, one difficulty seems to be that we
would need to analyze the expectation of moments or degree
distributions when the block model does not hold, which is doable. What
is worse is that the passage to $w$ from moments is very
ill-conditioned, involving first inversion via solution of the moment
problem, and then estimation of eigenvectors and eigenvalues from a
sequence of iterates $T_{w}(1),T_{w}^{2}(1)$, etc. If we assume
$\lambda_n\rightarrow\infty$ so that we can use consistency of the
degree distributions, we bypass the moment problem, but the
eigenfunction estimation problem remains. A step in this direction is a
result of \citet{Rohe2011} which shows that spectral clustering
can be
used to estimate the parameters of $k$ block models if $\lambda
\rightarrow\infty$ sufficiently, even if $k\rightarrow\infty$ slowly.
Unfortunately this does not deal with the problem we have just
discussed, how to pick a block model which is a good approximation to
the nonparametric model. For reasons which will appear in a future
paper, smoothness assumptions on $w$ have to be treated with caution.

While $\lambda_n \rightarrow\infty$ has not occurred in practice in
the past, networks with high average degrees are now
appearing routinely. In particular, university Facebook networks have
$\lambda$ of 15 or more with $n$ in the low
thousands. In any case $\lambda_n \rightarrow\infty$ can still be
useful as an asymptotic regime that can help us
understand some general patterns, in the same way that the sample size
going to infinity does in ordinary statistics.
Note that most of the time we do not specify the rate of growth of
$\lambda_n$, which can be very slow.

%s6.2 ###
\subsection{Adding covariates and directed graphs}

In principle, adding covariates $X_{i}$ at each vertex or $X_{ij}$
at each edge simply converts our latent variable model, $w(\cdot,\cdot)$
into a mixed model
\[
\mathbb{P}_{\theta}(A_{ij}=1|X_{i},X_j,X_{ij},\xi_{i},\xi_{j}) =
w_{\theta}(\xi_{i},\xi_{j},X_{i},X_j,X_{ij}),
\]
which can be turned into a logistic mixed model. Special cases of
such models have been considered in the literature; see \citet
{Hoff2007} and references therein. We do not pursue this
here. The extension of this model to directed graphs is also straightforward.

%s6.3 ###
\subsection{Dynamic models}

Many models
in the literature have been specified
dynamically; see \citet{Newman2010}. For instance, the
``preferential attachment''
model constructs an $n$ graph by
adding 1 vertex at a time, with edges of that vertex to previous vertices
formed with probabilities which are functions of the degree of
the candidate ``old'' vertex. If we let $n\rightarrow\infty$,
we obtain models of the type we have considered whose $w$ function
can be based on an integral equation for $\tau(\xi)$, our proxy for
the degree of the vertex with latent variable $\xi$. We shall pursue
this elsewhere also.

%We conjecture that for smooth functions $g$, \[
%n^{\frac{1}{2}}[\int g(\V{u})d\hat{F}_{m}(\V{u})-\int g(\V{u})d
% where $F_{m}^{*}=F_{m}$ if $\lambda_n\to\infty$, $F_{m}(\V{u})=\lim
% _{m}\mathbb{E}\hat{F}_{m}(\V{u})$
%if $\lambda_n=O(1)$.
%{\bf(shall we remove this or move to an earlier place?)}
%

\begin{appendix}\label{app}
\section*{Appendix: Additional lemmas and proofs}

\begin{pf*}{Proof of Proposition \ref{prop1}}
The first line of (\ref{eqT2}) is immediate, conditioning on $\{\xi
_{1},\ldots,\xi_{n}\}$. The second line in (\ref{eqT2}) follows by
expanding the second product. Finally, (\ref{eqT2}) follows directly
from the definitions of $P$ and~$Q$.
\end{pf*}

The following standard result is used in the proof of Theorem \ref{thmLandP}.
\begin{lem}
\label{lemcondindept} Suppose $(U_{n},V_{n})$ are random elements
such that,
\begin{eqnarray*}
\mathcal{L}(U_{n}) &\longrightarrow& \mathcal{L}(U),\\
\mathcal{L}(V_{n}|U_{n}) &\longrightarrow& \mathcal{L}(V)
\end{eqnarray*}
in probability. Then $U_{n}$, $V_{n}$ are asymptotically independent,
\[
\mathcal{L}(V_{n})\longrightarrow\mathcal{L}(V).
\]
\end{lem}
\begin{pf*}{Proof of Theorem \ref{thmLandP}}
By definition, $\mathbb{E} (\frac{L}{n\lambda_{n}})=\frac
{1}{2}$. Moreover,
\begin{eqnarray*}
\operatorname{Var} \biggl(\frac{1}{n\lambda_{n}} \sum\{A_{ij}\mbox{: all
}1\leq
i<j\leq n\} \biggr) &=&
(n\lambda_{n})^{-2} \mathbb{E} \biggl(\operatorname{Var}\biggl(\sum_{i < j}
A_{ij}\arrowvert\bolds{\xi}\biggr)\biggr)\\
&&{} +\rho_{n}^{2}(n\lambda_{n})^{-2}\operatorname{Var} \biggl( \sum_{i<j}w(\xi
_{i},\xi_{j})\biggr)\\
&\equiv&\operatorname{Var}(T_{1})+\operatorname{Var}(T_{2}) ,
\end{eqnarray*}
where
\begin{eqnarray*}
T_{1} & = & (n\lambda_{n})^{-1}\sum_{i<j}\bigl(A_{ij}-\rho_{n}w(\xi
_{i},\xi
_{j})\bigr) , \\
T_{2} & = & \rho_{n}(n\lambda_{n})^{-1}\sum_{i<j}w(\xi_{i},\xi
_{j})-\frac{1}{2} .
\end{eqnarray*}

Since $\lambda_{n}=(n-1)\rho_{n}$, the first term is
\begin{eqnarray*}
&&(n\lambda_{n})^{-2} \mathbb{E} \sum\bigl\{h(\xi_{i},\xi_{j})\bigl(1-h(\xi
_{i},\xi_{j})\bigr)\mbox{ all }i,j\bigr\}\\
&&\qquad\leq\frac{\rho
_{n}n^{2}}{2n^{2}\lambda
_{n}^{2}}
= O((n^{2}\rho_{n})^{-1})=O((n\lambda_{n})^{-1}) .
\end{eqnarray*}

The second term is a $U$-statistic of order 2, which is well known
to be $O(n^{-1})$. Thus, (\ref{eqneweq5}) follows in case (a).

To establish (\ref{eqneweq6}) and (b), we note that the conditional
distribution of $\sqrt{n\lambda_n}T_{1}$ given $\bolds{\xi}$ is that
of a sum of independent random variables with conditional variance
\[
\frac{1}{n\lambda_n}\sum_{i<j}\rho_{n}w(\xi_{i},\xi_{j})\bigl(1-\rho
_{n}w_{n}(\xi_{i},\xi_{j})\bigr) = \frac{1}{n^{2}}\sum_{i<j}w(\xi
_{i},\xi
_{j})\bigl(1+o_{P}(1)\bigr)
\stackrel{P}{\rightarrow} \frac{1}{2} .
\]
This sum is approximated by a $U$-statistic of order
2. Note that $\mathbb{E}w(\xi_{i},\xi_{j})=1$.
Since
the max of the summands in $\sqrt{n\lambda_n}T_{1}$ is $\frac
{1}{\sqrt
{n}\lambda_n}\rightarrow0$,
by the Lindeberg--Feller theorem, the conditional distribution tends
to $\mathcal{N}(0,\frac{1}{2})$ in probability. We can similarly
apply the limit theorem for $U$-statistics [see \citet{serfling}]
to conclude
that
\[
\sqrt{n}T_{2} \Rightarrow\mathcal{N}(0,\operatorname{Var}(\tau(\xi))).
\]
Applying Lemma \ref{lemcondindept}, we see that if $\lambda_n=O(1)$,
(b) follows. On the other hand, if $\lambda_n\rightarrow\infty$,
$\sqrt{n}T_{1}$
is negligible, and the Gaussian limit is determined by $T_{2}$.

The proof of (\ref{eqneweq7}) and (\ref{eqneweq8}) is similar.
We shall decompose $\check{P}(R)$ as $U_{1}+U_{2}$ as we did $\frac
{L}{n\lambda_n}$.
If $\lambda_n\rightarrow\infty$, it is enough to prove that
\[
\sqrt{n}\bigl(\check{P}(R)-\tilde{P}(R)\bigr)\Rightarrow\mathcal
{N}
(0,\sigma^{2}(R))
\]
since\vspace*{1pt} replacing $\bar{D}$ by $n\rho_{n}=\lambda_{n}$ gives
a perturbation of order $(n\lambda_{n})^{-{1/2}}=o(n^{-{1/2}})$.

In case (b), it is enough to show that the joint distribution of $\sqrt
{n}((\hat{P}(R)-P(R))\rho_{n}^{-|R|},T_{1},T_{2})$
is Gaussian\vspace*{1pt} in the limit, since in view of (\ref{eqneweq5}) and
(\ref{eqneweq6}) we can apply the delta method to $\check{P}(R)$.
Let $p\equiv|V(R)|$, $q\equiv|R|$. Each term in $\check{P}(R)$ is of
the form
\[
T(S)\equiv\frac{1}{{n \choose p}N(R)}\Pi\{
A_{i_{l}j_{l}}\dvtx(i_{l},j_{l})\in E(S), S\sim R\}.
\]
Condition on $\bolds{\xi} = \{\xi_{1},\ldots,\xi_{n}\}$. Then terms
$T(S)$, as
above, yield
%
%e6.1 ###
%
\begin{equation}\label{eq9}
\mathbb{E}(\hat{P}(R)|\bolds{\xi}) = \frac{1}{{n \choose
p}N(R)}\sum
_{S\sim R}\biggl(\prod_{(i,j)\in E(S)}[w(\xi_{i},\xi_{j})]\biggr) +
O(n^{-1}\lambda_n) .
\end{equation}
Thus,
\begin{eqnarray*}
U_{2} & = & \mathbb{E} (\hat{P}(R)|\bolds{\xi})\rho_{n}^{-q}-P(R),\\
U_{1} & = & \rho_{n}^{-q}\sum\{ T(S)-E(T(S)|\bolds{\xi}
)\dvtx S\sim
R\}.
\end{eqnarray*}

We begin by considering $\operatorname{Var}(U_{1}|\xi)$ which we can write as
\[
\sum\operatorname{cov}(T(S_{1}),T(S_{2})|\bolds{\xi})\rho_{n}^{-2q},
\]
where the sum ranges over all $S_{1}\sim R$, $S_{2}\sim R$.

If $E(S_{1})\cap E(S_{2})=\phi$ the covariance is 0. In general,
suppose the graph $S_{1}\cap S_{2}$ has $c$ vertices and $d$ edges.
Since $R$ is acyclic any subgraph is acyclic. By Corollary~3.2 of
\citet{chartrand86} for every acyclic graph,
$|V(S)| \geq|E(S)|+1$. Now,
%
%e6.2 ###
%
\begin{equation}\label{eq11}
\rho_{n}^{-2q}\operatorname{cov}(T(S_{1}),T(S_{2})|\bolds{\xi}) \leq
n^{-2p}\rho_{n}^{-d}\prod_{(i,j)\in S_{1}\cup S_{2}}w_{n}(\xi
_{i},\xi
_{j})
\end{equation}
since, if $d\geq1$,
%
%e6.3 ###
%
\begin{eqnarray}\label{eq12}
&&\mathbb{E} \Bigl[\Pi\{ A_{ij}\dvtx(i,j)\in\overline{S_{1}\cap
S_{2}}\}\Pi\{ A_{ij}^{2}\dvtx(i,j)\in S_{1}\cap S_{2}\}|\bolds
{\xi
}\Bigr]\nonumber\\[-8pt]\\[-8pt]
&&\qquad =\rho_{n}^{2q-d}\Pi\{ w_{n}(\xi_{i},\xi_{j})\dvtx(i,j)\in
S_{1}\cup
S_{2}\}.\nonumber
\end{eqnarray}
There are $O(n^{2p-c})$ terms in (\ref{eq9}) which have $c$ vertices
in common. Therefore by (\ref{eq11}) the total contribution of all
such terms to $\operatorname{Var}(U_{1})$ is
\[
O\biggl(n^{-c}\rho_{n}^{-d}\int w^{2q}(u,v)\,du \,dv\biggr),
\]
after using H\"older's inequality on
$\mathbb{E} \Pi\{ w(\xi_{i},\xi_{j})\dvtx(i,j)\in S_{1}\cup
S_{2}\}$.
From (\ref{eq12}) and our assumptions we conclude that
\[
\operatorname{Var}(U_{1})=O(n^{-1}\lambda_n^{-d})=o(n^{-1}),
\]
if $\lambda_n\to\infty$. On the other hand
\[
U_{2}=\frac{1}{{n \choose p}N(R)}\sum_{S\sim R}\biggl\{\prod_{(i,j)\in
S}w(\xi_{i},\xi_{j})\prod_{(i,j)\in\bar{S}}\bigl(1-h_{n}(\xi_{i},\xi
_{j})\bigr)-\tilde{P}(S)\biggr\}
\]
is a $U$-statistic. Its kernel
\[
\prod_{S}w(\xi_{i},\xi_{j})\prod_{\bar{S}}\bigl(1-h_{n}(\xi_{i},\xi
_{j})\bigr)-\tilde{P}(S)
\quad\stackrel{L_{2}}{\rightarrow}\quad \prod_{S}w(\xi_{i},\xi_{j})-\mathbb
{E}\prod_{S}w(\xi_{i},\xi_{j}).
\]
Thus, $\sqrt{n}(U_{1},U_{2})$ are jointly asymptotically Gaussian; see,
for instance, \citet{serfling}.

Since if $\lambda_n\to\infty$, $T_{1},U_{1}=o_{P}(n^{-{1/2}})$,
the result follows if $\lambda_n\to\infty$. If $\lambda_n=O(1)$, we
note that $\sqrt{n}(T_{1},U_{1})$ are sums of $q$ dependent random
variables in the sense of Bulinski [see \citet{doukhan94}] and hence,
given $\bolds{\xi}$, are jointly asymptotically Gaussian. It is not hard
to see that the limiting conditional covariance matrix is independent
of $\xi$, as it was for $T_{1}$ marginally. By Lemma \ref
{lemcondindept} again $(T_{1},U_{1})$
and $(T_{2},U_{2})$ are asymptotically independent and (a) and (b)
follow.\vspace*{1pt}

Finally\vspace*{1pt} we prove (c). To have $n^{-1/2}$ consistency for $\check{P}(R)$,
$\tilde{P}(R)$ and hence for $\check{Q}(R)$, $\tilde{Q}(R)$ by (\ref
{eqT4}) we need to argue that if $S\subset R$, $c\equiv|S|\leq p$
$|E(S)|=d$, then for a universal $M$,
\[
n^{-c}\rho^{-d}\leq Mn^{-1}.
\]
Since $\rho=\frac{\lambda_n}{n}$ we obtain
\[
n^{c} \biggl(\frac{\lambda_n}{n}\biggr)^{d} \geq n,\qquad
\lambda_n \geq n^{1-{(c-1)}/{d}}.
\]
For fixed $c\geq1$ this is maximized by $d=\frac{c(c-1)}{2}$ and
$n^{1-{2/c}}$ is maximized for $c\leq p$ by $c=p$.
\end{pf*}

\begin{pf*}{Proof of Theorem \ref{thm2}}
Since $T$ corresponds to the canonical $h$,
\begin{eqnarray*}
T(1)(\xi) & = & v_{(1)},\qquad 0\leq\xi\leq\pi_{1},\\
T(1)(\xi) & = & v_{(j)},\qquad \sum_{k=1}^{j-1}\pi_{k}\leq\xi\leq\sum
_{k=1}^{j}\pi_{k},\qquad 1\leq j\leq K,
\end{eqnarray*}
where $v_{(1)}<\cdots<v_{(k)}$ are the ordered $\{v_{j}\}$,
$v_{j}=\sum
_{i=1}^{K}\pi_{i}F_{ij}$.
By a theorem of Hausdorff and Hamburger [\citet{Feller2}], the
distribution
of the random variable $T(1)(\xi_{1})$ which takes on only $K$ distinct
values above is completely determined and uniquely so by its first
$2K-1$ moments $\mathbb{E} (T(1)(\xi_{1}))^{l}$, $l=1,\ldots,2K-1$. Therefore
for our model $\pi_{1},\ldots,\pi_{K}$ are completely determined
since $T(1)(\xi_{1})$ takes values $v_{j}$ with probability $\pi_{j}$,
$j=1,\ldots,K$.

Let $v^{(1)}=(v_{(1)},\ldots,v_{(K)})^{T}=F\pi$. Note that $\mathbb{E}
(T^{2}(1)(\xi_{1}))^{l},l=1,\ldots,2K-1$,
similarly determines the distribution of $T^{2}(1)(\xi_{1})$. Hence,
\[
v^{(2)} = Fv^{(1)}.
\]
Continuing we see that the $(K-1)(2K-1)$ moments $\{\tau_{kl}\dvtx2\leq
k\leq K,1\leq l\leq2K-1\}$
yield
%
%e6.4 ###
%
\begin{equation}\label{eqT7}
v^{(j)} = Fv^{(j-1)}
\end{equation}
for $j=1,\ldots,K$ where $v^{(0)}\equiv\pi$.

Given $\pi,v^{(1)},\ldots,v^{(K)}$ linearly independent, we can compute
$F$ since by (\ref{eqT7}), we can write
\[
F_{K\times K}V_{K\times K}^{(1)} = V_{K\times K}^{(2)},
\]
where $V^{(1)}=(v^{(0)},\ldots,v^{(K-1)})^{T}$ and
$V^{(2)}=(v^{(1)},\ldots,v^{(K)})^{T}$
and hence
\[
F = V^{(2)}\bigl[V^{(1)}\bigr]^{-1}.
\]
Consistency and $\sqrt{n}$-consistency follow from Theorem \ref{thmLandP}
and the delta method.
%\rightqed
\end{pf*}
\begin{pf*}{Proof of Proposition \ref{propA}}
Note that
%
%e6.5 ###
%
\begin{eqnarray}\label{eq17}
\mathbb{E} \exp sT^{l}(1)(\xi) & = & \mathbb{E} \exp s \mathbb{E}
\bigl(w(\xi,\xi_{1})\cdots
w(\xi_{l-1},\xi_{l})|\xi\bigr)\nonumber\\[-8pt]\\[-8pt]
& \leq & \mathbb{E} \exp s\bigl(w(\xi,\xi_{1})\cdots w(\xi_{l-1},\xi
_{l})\bigr).\nonumber
\end{eqnarray}
Taking $\xi=\xi_{0}$,
%
%e6.6 ###
%
\begin{equation}\label{eq18}
\mbox{(\ref{eq17})} \leq\mathbb{E}\exp|s|\Biggl(\frac{1}{l}\sum
_{j=0}^{l}w^{l}(\xi_{j},\xi_{j+1})\Biggr)
\end{equation}
by the arithmetic/geometric mean and Minkowski inequalities. By H\"older's
inequality (\ref{eq18}) is bounded by
\[
\prod_{j=0}^{l}[ \mathbb{E} \exp|s|w^{l}(\xi_{j},\xi
_{j+1})
]^{{1}/{l}}.
\]
It is easy to show that (A$'$) implies that $\mathbb{E} \exp\{
\sum
_{j=1}^{m}s_{j}T^{j}(1)(\xi)\}$
converges for $0<|s|<\varepsilon$ for some $\varepsilon$ depending on
$m$ and hence by a classical result that (A$'$) implies (A).
\end{pf*}
\begin{pf*}{Proof of Theorem \ref{thmnonparametric}}
Clearly $w$ determines the joint distribution of moments. We can take
$\tau_w(\xi) = T_w(1)(\xi)$ monotone, corresponding to the
canonical~$w$, to be the quantile function of the marginal distribution
of $T_w(1)(\xi)$. Now the joint distribution of $(T_w(1)(\xi),
T_w^2(1)(\xi))$ determines $\tau_w(\cdot)$, $T_w \tau_w(\cdot)$, except
on a set of measure 0. Continuing this argument, we can determine the
entire sequence of functions $\tau _w$, $T_w\tau_w$, $T_w^2\tau_w,
\ldots.$ Since $T_w$ is bounded self-adjoint, these functions are all
in $L_2$. Let $g_{k}^{(1)}(\cdot)=T_w(\frac
{g_{k-1}^{(1)}}{|g_{k-1}^{(1)}|})$, $g_{0}^{(1)}(\cdot)=1$, where $|f|$
and $(f,g)$ are, respectively, the norm and the inner product in $L_2$.
Then $g_{k}\rightarrow_{L_{2}}\lambda_{1}\phi_{1}$ where $\lambda_{1}$
is the first eigenvalue, $\phi_{1}$ the first eigenfunction and
$\frac{g_{k}}{|g_{k}|}\rightarrow\phi_{1}$. This is just the ``powering
up'' method applied to the function 1 with convergence guaranteed since
$\lambda_1$ is unique, and 1 is not orthogonal to $\phi_1$ or any other
eigenfunction. So $\lambda_{1}$ and $\phi_{1}$ are also determined.
Thus we can compute $g_{0}^{(2)}\equiv1-(1,\phi_{1})\phi_{1}$. Further,
\[
g_{1}^{(2)}=T_w\biggl(\frac{g_{0}^{(2)}}{|g_{0}^{(2)}|}\biggr)=
\frac{T_w1(\cdot)-\lambda_{1}(1,\phi_{1})\phi_1}{|1-(1,\phi
_{1})\phi_{1}|}
\]
is computable since we know $T_w1(\cdot)$ and the eigenfunction $\phi_{1}$
and eigenvalue $\lambda_{1}$. More generally, $T_w^{k}g_{1}^{(2)}$,
$|g_{k-1}^{(2)}|$ can be similarly determined. Then, by the same
argument as before, using 1 not orthogonal to $\phi_2$, we obtain
$g_{k}^{(1)}\rightarrow_{L_{2}}\lambda_{2}\phi_{2}$
and $g_{k}^{(1)}/|g_{k}^{(1)}|\rightarrow_{L_{2}}\phi_{2}$. Now form
$g_{0}^{(3)}\equiv1-\lambda_{1}(1,\phi_{1})\phi_{1}-\lambda
_{2}(1,\phi
_{2})\phi_{2}$
and proceed as before, and continue to determine $\lambda_k, \phi_k$
for all $k$. This and (\ref{eigenw}) complete the proof.
\end{pf*}
\begin{pf*}{Proof of Theorem \ref{thmldegree}}
Note first that (\ref{eqmallowsjoint}) implies that the $M_2$ distance
between $\hat F_m$ and the empirical distribution of $\{ \bolds{
\theta}_m(\xi_i) \}$ tends to 0. The first conclusion of the theorem
now follows by the Glivenko--Cantelli theorem and the Law of Large Numbers.

To show (\ref{eqmallowsjoint}), note that
%
%e6.7 ###
%
\begin{equation}\label{eq19}
\frac{1}{n}\sum_{i=1}^{n}\bigl|\tilde{D}_{i}^{(m)}-\theta_{m}(\xi_{i})\bigr|^{2}
\stackrel{P}{\rightarrow} 0,
\end{equation}
where $\tilde{D}_{i}^{(m)}\equiv(\frac{D_{i}}{\bar{D}},\ldots
,\frac
{D_{i}^{(m)}}{\bar{D}^{m}})^{T}$.
By Theorem \ref{thmLandP}, we can replace $\bar{D}$ by $\lambda_n$ if
$\lambda_n\geq\varepsilon$.
Then (\ref{eq19}) is implied by
%
%e6.8 ###
%
\begin{equation}\label{eq20}
\frac{1}{n}\sum_{i=1}^{n}\mathbb{E}\Biggl|
\sum_{j=1}^{n}\frac{\tilde{A}_{ij}^{(m)}}{\lambda_n^{m}}-\theta
_{m}(\xi_{i}) \Biggr|^{2} \rightarrow0.
\end{equation}
Now,
%
%e6.9 ###
%
\begin{eqnarray}\label{eq21}
&&\sum_{j=1}^{n}\mathbb{E}\biggl(\frac{\tilde{A}_{ij}^{(m)}}{\lambda
_n^{m}} \Big|
\bolds{\xi}\biggr) \nonumber\\
&&\qquad = \frac{1}{n^{m}}\sum\bigl\{w_{E(R)}\dvtx R=\{
(i,i_{1}),\ldots,(i_{m-1},j)\},\\
&&\hspace*{140pt}\mbox{all vertices distinct}\bigr\},\nonumber
\end{eqnarray}
where $w_{E(R)} = \prod_{(a,b)\in E(R)}w(\xi_{a},\xi_{b})$. Further,
(\ref{eq21}) is a $U$-statistic of order $m$ under $|w_{2m}|<\infty$
and
\[
\mathbb{E}\Biggl|\sum_{j=1}^{n}\mathbb{E}\biggl(\frac{\tilde
{A}_{ij}^{(m)}}{\lambda_n^{m}}\Big|
\bolds{\xi}\biggr)-\mathbb{E}\bigl(w_{E(R)}|\xi_{i}\bigr)\Biggr|^{2} \leq\frac
{C|w_{2m}|}{n}
\]
by standard theory [\citet{serfling}].

Since $\mathbb{E}(w_{E(R)}|\xi_{i})=\theta_{m}(\xi_{i})$, we can
consider
%
%e6.10 ###
%
\begin{eqnarray}\label{eq22}
&&\mathbb{E} \Biggl(\frac{1}{n}\sum_{i=1}^{n}\Biggl|\sum
_{j=1}^{n}\frac
{\tilde{A}_{ij}^{(m)}-\mathbb{E}(\tilde{A}_{ij}^{(m)}|
\bolds{\xi})}{\lambda_n^{m}}\Biggr|^{2}\Biggr) \nonumber\\[-8pt]\\[-8pt]
&&\qquad\leq\max_{i} \frac{\mathbb{E}|\sum_{j=1}^{n} (\tilde
{A}_{ij}^{(m)}-\mathbb{E}(\tilde{A}_{ij}^{(m)}|
\bolds{\xi}))|^{2}}{\lambda_n^{2m}}.\nonumber
\end{eqnarray}
Note that $R = \{(i,i_{1}),(i_{1},i_{2}),\ldots,(i_{m-1},j)\}$
is acyclic if all vertices are distinct.
As in the proof of Theorem \ref{thmLandP},
all nonzero covariance terms in (\ref{eq22}) are of order $\rho
^{2m-d}n^{2m-c}$
where $c\geq d$ since the intersection graphs all have $i$ in common
but are otherwise acyclic. The largest order term corresponds to $c=d=m$,
so that
\[
\mathbb{E}\Biggl|\sum_{j=1}^{n}\bigl(\lambda_n^{-m}\tilde
{A}_{ij}^{(m)}-\theta_{m}(\xi_{i})\bigr)\Biggr|^{2} \leq C\lambda_n^{-m},
\]
where $C$ depends on $|w_{2m}|$ only. Thus (\ref{eq20}) holds
if $\lambda_n\rightarrow\infty$.
\end{pf*}
\end{appendix}

\section*{Acknowledgment}

Thanks to Allan Sly for a helpful discussion.

%suskaldyti doi

% imsref loaded by lrinkeviciute, 2011-09-13 16:10:35
% imsref loaded by lrinkeviciute, 2011-09-13 16:16:56
%
% imsref loaded by lrinkeviciute, 2011-10-05 13:24:46

\printaddresses


\begin{thebibliography}{30}
% BibTex style file: ims.bst, 2011-05-30
% Default style options (sort=0,type=number).
% Used options (sort=1,type=nameyear).

%b1 ###
\bibitem[\protect\citeauthoryear{Aldous}{1981}]{Aldous81}
%
\begin{barticle}[mr]
\bauthor{\bsnm{Aldous},~\bfnm{David~J.}\binits{D.~J.}}
(\byear{1981}).
\btitle{Representations for partially exchangeable arrays of random variables}.
\bjournal{J.~Multivariate Anal.}
\bvolume{11}
\bpages{581--598}.
\bid{doi={10.1016/0047-259X(81)90099-3}, issn={0047-259X}, mr={0637937}}
\bptok{imsref}%
\end{barticle}
%
\endbibitem

%b2 ###
\bibitem[\protect\citeauthoryear{Barab{\'a}si and
Albert}{1999}]{BarabasiAlbert1999}
%
\begin{barticle}[mr]
\bauthor{\bsnm{Barab{\'a}si},~\bfnm{Albert-L{\'a}szl{\'o}}\binits
{A.-L.}} \AND
\bauthor{\bsnm{Albert},~\bfnm{R{\'e}ka}\binits{R.}}
(\byear{1999}).
\btitle{Emergence of scaling in random networks}.
\bjournal{Science}
\bvolume{286}
\bpages{509--512}.
\bid{doi={10.1126/science.286.5439.509}, issn={0036-8075}, mr={2091634}}
\bptok{imsref}%
\end{barticle}
%
\endbibitem

%b3 ###
\bibitem[\protect\citeauthoryear{Bickel and Chen}{2009}]{BickelChen2009}
%
\begin{barticle}[pbm]
\bauthor{\bsnm{Bickel},~\bfnm{Peter~J.}\binits{P.~J.}} \AND
\bauthor{\bsnm{Chen},~\bfnm{Aiyou}\binits{A.}}
(\byear{2009}).
\btitle{A nonparametric view of network models and Newman--Girvan and other
modularities}.
\bjournal{Proc. Natl. Acad. Sci. USA}
\bvolume{106}
\bpages{21068--21073}.
\bid{doi={10.1073/pnas.0907096106}, issn={1091-6490}, pii={0907096106},
pmcid={2795514}, pmid={19934050}}
\bptok{imsref}%
\end{barticle}
%
\endbibitem

%b4 ###
\bibitem[\protect\citeauthoryear{Bollob{\'a}s, Janson and
Riordan}{2007}]{Bollobas2007}
%
\begin{barticle}[mr]
\bauthor{\bsnm{Bollob{\'a}s},~\bfnm{B{\'e}la}\binits{B.}},
\bauthor{\bsnm{Janson},~\bfnm{Svante}\binits{S.}} \AND
\bauthor{\bsnm{Riordan},~\bfnm{Oliver}\binits{O.}}
(\byear{2007}).
\btitle{The phase transition in inhomogeneous random graphs}.
\bjournal{Random Structures Algorithms}
\bvolume{31}
\bpages{3--122}.
\bid{doi={10.1002/rsa.20168}, issn={1042-9832}, mr={2337396}}
\bptok{imsref}%
\end{barticle}
%
\endbibitem

%b5 ###
\bibitem[\protect\citeauthoryear{Chartrand, Lesniak and
Behzad}{1986}]{chartrand86}
%
\begin{bbook}[auto:STB|2011/09/12|07:03:23]
\bauthor{\bsnm{Chartrand},~\bfnm{G.}\binits{G.}},
\bauthor{\bsnm{Lesniak},~\bfnm{L.}\binits{L.}} \AND
\bauthor{\bsnm{Behzad},~\bfnm{M.}\binits{M.}}
(\byear{1986}).
\btitle{Graphs and Digraphs}, \bedition{2nd} ed.
\bpublisher{Wadsworth and Brooks}, \baddress{Monterey, CA}.
\bptok{imsref}%
\end{bbook}
%
\endbibitem

%b6 ###
\bibitem[\protect\citeauthoryear{Chatterjee and
Diaconis}{2011}]{chatterjeediaconis2011}
%
\begin{bmisc}[auto:STB|2011/09/12|07:03:23]
\bauthor{\bsnm{Chatterjee},~\bfnm{S.}\binits{S.}} \AND
\bauthor{\bsnm{Diaconis},~\bfnm{P.}\binits{P.}}
(\byear{2011}).
\bhowpublished{Estimating and understanding exponential
random graph models.
Unpublished manuscript.
Available at \href{http://arxiv.org/abs/arXiv:1102.2650}{arXiv:1102.2650}.}
\bptok{imsref}%
\end{bmisc}
%
\endbibitem

%b7 ###
\bibitem[\protect\citeauthoryear{Chung and Lu}{2002}]{ChungLu2002}
%
\begin{barticle}[mr]
\bauthor{\bsnm{Chung},~\bfnm{Fan}\binits{F.}} \AND
\bauthor{\bsnm{Lu},~\bfnm{Linyuan}\binits{L.}}
(\byear{2002}).
\btitle{Connected components in random graphs with given expected degree
sequences}.
\bjournal{Ann. Comb.}
\bvolume{6}
\bpages{125--145}.
\bid{doi={10.1007/PL00012580}, issn={0218-0006}, mr={1955514}}
\bptok{imsref}%
\end{barticle}
%
\endbibitem

%b8 ###
\bibitem[\protect\citeauthoryear{de~Solla~Price}{1965}]{price1965}
%
\begin{barticle}[auto:STB|2011/09/12|07:03:23]
\bauthor{\bparticle{de} \bsnm{Solla~Price},~\bfnm{D.~J.}\binits{D.~J.}}
(\byear{1965}).
\btitle{Networks of scientific papers}.
\bjournal{Science}
\bvolume{149}
\bpages{510--515}.
\bptok{imsref}%
\end{barticle}
%
\endbibitem

\bibitem[\protect\citeauthoryear{Decelle et al.}{2011}]{Decelleetal2011}
\begin{bmisc}[auto]
\bauthor{\bsnm{Decelle},~\bfnm{Aurelien}\binits{A.}},
\bauthor{\bsnm{Krzakala},~\bfnm{Florent}\binits{F.}},
\bauthor{\bsnm{Moore},~\bfnm{Cristopher}\binits{C.}} \AND
\bauthor{\bsnm{Zdeborov\'{a}},~\bfnm{Lenka}\binits{L.}}
(\byear{2011}).
\bhowpublished{Asymptotic analysis of the stochastic block
model for modular networks and its algorithmic applications.
Available at \href{http://arxiv.org/abs/1109.3041}{arXiv:1109.3041}.}
\bptok{imsref}%
\end{bmisc}
%
\endbibitem

%b9 ###
\bibitem[\protect\citeauthoryear{Diaconis and Janson}{2008}]{Diaconis2008}
%
\begin{barticle}[mr]
\bauthor{\bsnm{Diaconis},~\bfnm{Persi}\binits{P.}} \AND
\bauthor{\bsnm{Janson},~\bfnm{Svante}\binits{S.}}
(\byear{2008}).
\btitle{Graph limits and exchangeable random graphs}.
\bjournal{Rend. Mat. Appl. (7)}
\bvolume{28}
\bpages{33--61}.
\bid{issn={1120-7183}, mr={2463439}}
\bptok{imsref}%
\end{barticle}
%
\endbibitem

%b10 ###
\bibitem[\protect\citeauthoryear{Doukhan}{1994}]{doukhan94}
%
\begin{bbook}[mr]
\bauthor{\bsnm{Doukhan},~\bfnm{Paul}\binits{P.}}
(\byear{1994}).
\btitle{Mixing: Properties and Examples}.
\bseries{Lecture Notes in Statistics}
\bvolume{85}.
\bpublisher{Springer}, \baddress{New York}.
\bid{mr={1312160}}
\bptok{imsref}%
\end{bbook}
%
\endbibitem

%b11 ###
\bibitem[\protect\citeauthoryear{Feller}{1971}]{Feller2}
%
\begin{bbook}[mr]
\bauthor{\bsnm{Feller},~\bfnm{William}\binits{W.}}
(\byear{1971}).
\btitle{An Introduction to Probability Theory and Its Applications. {V}ol.
{II}},
\bedition{2nd} ed.
\bpublisher{Wiley}, \baddress{New York}.
\bid{mr={0270403}}
\bptok{imsref}%
\end{bbook}
%
\endbibitem

%b12 ###
\bibitem[\protect\citeauthoryear{Frank and Strauss}{1986}]{FrankStrauss1986}
%
\begin{barticle}[mr]
\bauthor{\bsnm{Frank},~\bfnm{Ove}\binits{O.}} \AND
\bauthor{\bsnm{Strauss},~\bfnm{David}\binits{D.}}
(\byear{1986}).
\btitle{Markov graphs}.
\bjournal{J. Amer. Statist. Assoc.}
\bvolume{81}
\bpages{832--842}.
\bid{issn={0162-1459}, mr={0860518}}
\bptok{imsref}%
\end{barticle}
%
\endbibitem

%b13 ###
\bibitem[\protect\citeauthoryear{Handcock}{2003}]{handcock2003}
%
\begin{bmisc}[auto:STB|2011/09/12|07:03:23]
\bauthor{\bsnm{Handcock},~\bfnm{M.}\binits{M.}}
(\byear{2003}).
\bhowpublished{Assessing degeneracy in statistical models of social
networks. Working Paper 39,
Center for Statistics and the Social Sciences.}
\bptok{imsref}%
\end{bmisc}
%
\endbibitem

%b14 ###
\bibitem[\protect\citeauthoryear{Handcock, Raftery and
Tantrum}{2007}]{Handcock2007}
%
\begin{barticle}[mr]
\bauthor{\bsnm{Handcock},~\bfnm{Mark~S.}\binits{M.~S.}},
\bauthor{\bsnm{Raftery},~\bfnm{Adrian~E.}\binits{A.~E.}} \AND
\bauthor{\bsnm{Tantrum},~\bfnm{Jeremy~M.}\binits{J.~M.}}
(\byear{2007}).
\btitle{Model-based clustering for social networks}.
\bjournal{J. Roy. Statist. Soc. Ser. A}
\bvolume{170}
\bpages{301--354}.
\bid{doi={10.1111/j.1467-985X.2007.00471.x}, issn={0964-1998}, mr={2364300}}
\bptok{imsref}%
\end{barticle}
%
\endbibitem

%b15 ###
\bibitem[\protect\citeauthoryear{Hoff}{2007}]{Hoff2007}
%
\begin{bincollection}[auto:STB|2011/09/12|07:03:23]
\bauthor{\bsnm{Hoff},~\bfnm{P.~D.}\binits{P.~D.}}
(\byear{2007}).
\btitle{Modeling homophily and stochastic equivalence in symmetric relational
data}.
In \bbooktitle{Advances in Neural Information Processing Systems}
\bvolume{19}.
\bpublisher{MIT Press}, \baddress{Cambridge, MA}.
\bptok{imsref}%
\end{bincollection}
%
\endbibitem

%b16 ###
\bibitem[\protect\citeauthoryear{Hoff, Raftery and Handcock}{2002}]{Hoff2002}
%
\begin{barticle}[mr]
\bauthor{\bsnm{Hoff},~\bfnm{Peter~D.}\binits{P.~D.}},
\bauthor{\bsnm{Raftery},~\bfnm{Adrian~E.}\binits{A.~E.}} \AND
\bauthor{\bsnm{Handcock},~\bfnm{Mark~S.}\binits{M.~S.}}
(\byear{2002}).
\btitle{Latent space approaches to social network analysis}.
\bjournal{J. Amer. Statist. Assoc.}
\bvolume{97}
\bpages{1090--1098}.
\bid{doi={10.1198/016214502388618906}, issn={0162-1459}, mr={1951262}}
\bptok{imsref}%
\end{barticle}
%
\endbibitem

%b17 ###
\bibitem[\protect\citeauthoryear{Holland, Laskey and
Leinhardt}{1983}]{Holland83}
%
\begin{barticle}[mr]
\bauthor{\bsnm{Holland},~\bfnm{Paul~W.}\binits{P.~W.}},
\bauthor{\bsnm{Laskey},~\bfnm{Kathryn~Blackmond}\binits{K.~B.}}
\AND
\bauthor{\bsnm{Leinhardt},~\bfnm{Samuel}\binits{S.}}
(\byear{1983}).
\btitle{Stochastic blockmodels: First steps}.
\bjournal{Social Networks}
\bvolume{5}
\bpages{109--137}.
\bid{doi={10.1016/0378-8733(83)90021-7}, issn={0378-8733}, mr={0718088}}
\bptok{imsref}%
\end{barticle}
%
\endbibitem

%b18 ###
\bibitem[\protect\citeauthoryear{Holland and
Leinhardt}{1981}]{HollandLeinhardt1981}
%
\begin{barticle}[mr]
\bauthor{\bsnm{Holland},~\bfnm{Paul~W.}\binits{P.~W.}} \AND
\bauthor{\bsnm{Leinhardt},~\bfnm{Samuel}\binits{S.}}
(\byear{1981}).
\btitle{An exponential family of probability distributions for directed
graphs}.
\bjournal{J. Amer. Statist. Assoc.}
\bvolume{76}
\bpages{33--65}.
%Wasserman, Ove Frank and Shelby J. Haberman and a reply by the authors}.
\bid{issn={0162-1459}, mr={0608176}}
\bptnote{check related}%
\bptok{imsref}%
\end{barticle}
%
\endbibitem

%b19 ###
\bibitem[\protect\citeauthoryear{Hoover}{1979}]{Hoover1979}
%
\begin{bmisc}[auto:STB|2011/09/12|07:03:23]
\bauthor{\bsnm{Hoover},~\bfnm{D.}\binits{D.}}
(\byear{1979}).
\bhowpublished{Relations on probability spaces and arrays of random variables.
Technical report, Institute for Advanced Study, Princeton, NJ}.
\bptok{imsref}%
\end{bmisc}
%
\endbibitem

%b20 ###
\bibitem[\protect\citeauthoryear{Kallenberg}{2005}]{Kallenberg2005}
%
\begin{bbook}[mr]
\bauthor{\bsnm{Kallenberg},~\bfnm{Olav}\binits{O.}}
(\byear{2005}).
\btitle{Probabilistic Symmetries and Invariance Principles}.
\bpublisher{Springer}, \baddress{New York}.
\bid{mr={2161313}}
\bptok{imsref}%
\end{bbook}
%
\endbibitem

%b21 ###
\bibitem[\protect\citeauthoryear{Karrer and Newman}{2011}]{Karrer10}
%
\begin{barticle}[mr]
\bauthor{\bsnm{Karrer},~\bfnm{Brian}\binits{B.}} \AND
\bauthor{\bsnm{Newman},~\bfnm{M.~E.~J.}\binits{M.~E.~J.}}
(\byear{2011}).
\btitle{Stochastic blockmodels and community structure in networks}.
\bjournal{Phys. Rev. E (3)}
\bvolume{83}
\bpages{016107}.
\bid{doi={10.1103/PhysRevE.83.016107}, issn={1539-3755}, mr={2788206}}
\bptok{imsref}%
\end{barticle}
%
\endbibitem

%b22 ###
\bibitem[\protect\citeauthoryear{Newman}{2006}]{Newman2006}
%
\begin{barticle}[mr]
\bauthor{\bsnm{Newman},~\bfnm{M.~E.~J.}\binits{M.~E.~J.}}
(\byear{2006}).
\btitle{Finding community structure in networks using the eigenvectors of
matrices}.
\bjournal{Phys. Rev. E (3)}
\bvolume{74}
\bpages{036104}.
\bid{doi={10.1103/PhysRevE.74.036104}, issn={1539-3755}, mr={2282139}}
\bptok{imsref}%
\end{barticle}
%
\endbibitem

%b23 ###
\bibitem[\protect\citeauthoryear{Newman}{2010}]{Newman2010}
%
\begin{bbook}[mr]
\bauthor{\bsnm{Newman},~\bfnm{M.~E.~J.}\binits{M.~E.~J.}}
(\byear{2010}).
\btitle{Networks: An Introduction}.
\bpublisher{Oxford Univ. Press}, \baddress{Oxford}.
\bid{doi={10.1093/acprof:oso/9780199206650.001.0001}, mr={2676073}}
\bptok{imsref}%
\end{bbook}
%
\endbibitem

%b24 ###
\bibitem[\protect\citeauthoryear{Nowicki and Snijders}{2001}]{Nowicki2001}
%
\begin{barticle}[mr]
\bauthor{\bsnm{Nowicki},~\bfnm{Krzysztof}\binits{K.}} \AND
\bauthor{\bsnm{Snijders},~\bfnm{Tom A.~B.}\binits{T.~A.~B.}}
(\byear{2001}).
\btitle{Estimation and prediction for stochastic blockstructures}.
\bjournal{J. Amer. Statist. Assoc.}
\bvolume{96}
\bpages{1077--1087}.
\bid{doi={10.1198/016214501753208735}, issn={0162-1459}, mr={1947255}}
\bptok{imsref}%
\end{barticle}
%
\endbibitem

%b25 ###
\bibitem[\protect\citeauthoryear{Picard et~al.}{2008}]{Picardetal2008}
\begin{barticle}[mr]
\bauthor{\bsnm{Picard},~\bfnm{F.}\binits{F.}},
  \bauthor{\bsnm{Daudin},~\bfnm{J.~J.}\binits{J.~J.}},
  \bauthor{\bsnm{Koskas},~\bfnm{M.}\binits{M.}},
  \bauthor{\bsnm{Schbath},~\bfnm{S.}\binits{S.}} \AND
  \bauthor{\bsnm{Robin},~\bfnm{S.}\binits{S.}}
(\byear{2008}).
\btitle{Assessing the exceptionality of network motifs}.
\bjournal{J. Comput. Biol.}
\bvolume{15}
\bpages{1--20}.
\bid{issn={1066-5277}, mr={2383618}}
\bptok{imsref}%
\end{barticle}
\endbibitem

%b26 ###
\bibitem[\protect\citeauthoryear{Robins et~al.}{2007}]{robinsetc2007}
%
\begin{barticle}[auto:STB|2011/09/12|07:03:23]
\bauthor{\bsnm{Robins},~\bfnm{G.}\binits{G.}},
\bauthor{\bsnm{Snijders},~\bfnm{T.}\binits{T.}},
\bauthor{\bsnm{Wang},~\bfnm{P.}\binits{P.}},
\bauthor{\bsnm{Handcock},~\bfnm{M.}\binits{M.}} \AND
\bauthor{\bsnm{Pattison},~\bfnm{P.}\binits{P.}}
(\byear{2007}).
\btitle{Recent developments in exponential random graphs models
($p^*$) for
social networks}.
\bjournal{Social Networks}
\bvolume{29}
\bpages{192--215}.
\bptok{imsref}%
\end{barticle}
%
\endbibitem

%b27 ###
\bibitem[\protect\citeauthoryear{Rohe, Chatterjee and Yu}{2011}]{Rohe2011}
%
\begin{bmisc}[auto:STB|2011/09/12|07:03:23]
\bauthor{\bsnm{Rohe},~\bfnm{K.}\binits{K.}},
\bauthor{\bsnm{Chatterjee},~\bfnm{S.}\binits{S.}} \AND
\bauthor{\bsnm{Yu},~\bfnm{B.}\binits{B.}}
(\byear{2011}).
\bhowpublished{Spectral clustering and the high-dimensional stochastic block model.
\textit{Ann. Statist.} To appear.}
\bptok{imsref}%
\end{bmisc}
%
\endbibitem

%b28 ###
\bibitem[\protect\citeauthoryear{Serfling}{1980}]{serfling}
%
\begin{bbook}[mr]
\bauthor{\bsnm{Serfling},~\bfnm{Robert~J.}\binits{R.~J.}}
(\byear{1980}).
\btitle{Approximation Theorems of Mathematical Statistics}.
\bpublisher{Wiley}, \baddress{New York}.
\bid{mr={0595165}}
\bptok{imsref}%
\end{bbook}
%
\endbibitem

%b29 ###
\bibitem[\protect\citeauthoryear{Shalizi and
Rinaldo}{2011}]{ShaliziRinaldo2011}
%
\begin{bmisc}[auto:STB|2011/09/12|07:03:23]
\bauthor{\bsnm{Shalizi},~\bfnm{C.~R.}\binits{C.~R.}} \AND
\bauthor{\bsnm{Rinaldo},~\bfnm{A.}\binits{A.}}
(\byear{2011}).
\bhowpublished{Projective structure and parametric inference in exponential
families. Carnegie Mellon Univ. Unpublished manuscript.}
\bptok{imsref}%
\end{bmisc}
%
\endbibitem

%b30 ###
\bibitem[\protect\citeauthoryear{Snijders and
Nowicki}{1997}]{SnijdersNowicki1997}
%
\begin{barticle}[mr]
\bauthor{\bsnm{Snijders},~\bfnm{Tom A.~B.}\binits{T.~A.~B.}} \AND
\bauthor{\bsnm{Nowicki},~\bfnm{Krzysztof}\binits{K.}}
(\byear{1997}).
\btitle{Estimation and prediction for stochastic blockmodels for
graphs with
latent block structure}.
\bjournal{J. Classification}
\bvolume{14}
\bpages{75--100}.
\bid{doi={10.1007/s003579900004}, issn={0176-4268}, mr={1449742}}
\bptok{imsref}%
\end{barticle}
%
\endbibitem

\end{thebibliography}
\end{document}